\documentclass[twocolumn]{article}
\usepackage{framed,multirow}
\usepackage{enumerate}
\usepackage[colorlinks=true,linkcolor=blue,citecolor=blue,urlcolor=blue]{hyperref}
\usepackage{amssymb}
\usepackage{latexsym}
\usepackage{url}
\usepackage{xcolor}
\definecolor{newcolor}{rgb}{.8,.349,.1}
\definecolor{lightblue}{rgb}{0.6, 0.81, 0.93}  
\usepackage{subcaption}
\usepackage{authblk} 
\usepackage{graphicx}
\usepackage{svg}
\usepackage{booktabs}
\usepackage{amsmath}
\DeclareMathOperator*{\argmax}{argmax}
\DeclareMathOperator*{\argmin}{argmin}
\usepackage{pifont}
\newcommand{\cmark}{\ding{51}}

\newcommand{\rv}[1]{\textcolor{black}{#1}}

\title{Meta-Learning from Learning Curves for Budget-Limited Algorithm Selection \\
\vspace{0.5em}
\large \textit{Manuscript published in Pattern Recognition Letters, September 2024}}

\author[1]{Manh Hung Nguyen\thanks{Corresponding author: manh.hung.nguyen@chalearn.org}}
\author[2]{Lisheng Sun-Hosoya}
\author[1,2]{Isabelle Guyon}

\affil[1]{Chalearn, California, USA}
\affil[2]{TAU Team, LISN, Université Paris-Saclay, Gif-sur-Yvette, France}

\date{} 
\begin{document}
\maketitle

\begin{abstract}
Training a large set of machine learning algorithms to convergence in order to select the best-performing algorithm for a dataset is computationally wasteful. Moreover, in a budget-limited scenario, it is crucial to carefully select an algorithm candidate and allocate a budget for training it, ensuring that the limited budget is optimally distributed to favor the most promising candidates. Casting this problem as a Markov Decision Process, we propose a novel framework in which an agent must select in the process of learning the most promising algorithm without waiting until it is fully trained. At each time step, given an observation of partial learning curves of algorithms, the agent must decide whether to allocate resources to further train the most promising algorithm (exploitation), to wake up another algorithm previously put to sleep, or to start training a new algorithm (exploration). In addition, our framework allows the agent to meta-learn from learning curves on past datasets along with dataset meta-features and algorithm hyperparameters. By incorporating meta-learning, we aim to avoid myopic decisions based solely on premature learning curves on the dataset at hand. We introduce two benchmarks of learning curves that served in international competitions at WCCI'22 and AutoML-conf'22, of which we analyze the results. Our findings show that both meta-learning and the progression of learning curves enhance the algorithm selection process, as evidenced by methods of winning teams and our DDQN baseline, compared to heuristic baselines or a random search. Interestingly, our cost-effective baseline, which selects the best-performing algorithm w.r.t. a small budget, can perform decently when learning curves do not intersect frequently.
\end{abstract}



\section{Introduction}
In a typical Machine Learning (ML) task, given a dataset, one is asked to build a model for the dataset w.r.t. an objective (e.g., classification). The model is often defined as the output of training an ML algorithm (e.g., neural network weights) on the dataset, capturing the learned data distribution. However, selecting the best-suited algorithm for a particular dataset is challenging, especially for non-experts with limited ML knowledge. This \textit{algorithm selection} problem becomes even more challenging in a \textit{budget-limited} scenario, where evaluating a large set of algorithms becomes costly if each of them must be trained or optimized to convergence. As a consequence, a common practice in the ML community is early stopping or discarding unpromising algorithms based on their \textit{learning curves} during training [\cite{DBLP:journals/corr/abs-2201-12150}]. A learning curve evaluates an algorithm's incremental performance improvement, as a function of time, number of epochs, or number of training examples. Learning how to effectively leverage such information to boost algorithm selection is the focus of this paper.

{In the literature, a common approach is to invest time in collecting premature learning curves of algorithm candidates and extrapolate them to select the most promising algorithm [\cite{DBLP:conf/ijcai/DomhanSH15, DBLP:journals/corr/abs-2201-12150, DBLP:conf/ijcnn/SchmidtGNS20, DBLP:journals/pami/MohrR23}]. Notably, several existing works use neural networks for this purpose. For example, \cite{DBLP:conf/iclr/KleinFSH17} use Bayesian neural networks for modeling and predicting learning curves, while \cite{DBLP:conf/nips/AdriaensenR0H23} use prior-data fitted networks trained to extrapolate artificial right-censored learning curves generated from a parametric prior. However, these approaches often require Parametric Learning Curve Models and some of them rely on the assumption of learning curve concavity. Moreover, they suffer from myopia due to extrapolations based solely on partially observed curves without considering meta-knowledge. }

{\textit{Meta-learning} can be incorporated to address these issues. Meta-learning has shown great potential in learning from previous tasks to solve new ones more efficiently [\cite{DBLP:journals/corr/abs-1810-03548, vettoruzzo2024advances}]. In scenarios where meta-knowledge, such as features of datasets or past performances of algorithms, is available, meta-learning can help improve the algorithm selection process [\cite{DBLP:journals/ml/BrazdilG18, Brazdil2022, DBLP:journals/jmlr/FeurerEFLH22}]. For instance, \cite{DBLP:conf/pkdd/AbdulrahmanBRV15} uses meta-level information acquired in past experiments to construct an average ranking of algorithms and apply active testing. \cite{DBLP:conf/ida/RijnABV15} exploit the similarity of the partially-observed rankings, and use the most similar learning curve as surrogates. \cite{DBLP:conf/icml/WistubaP20} propose a method that learns to rank learning curves by optimizing a pairwise ranking loss. More recent works use LSTM networks and transformer-based models to learn embeddings of dataset meta-features and algorithm candidates' learning behavior observed on other datasets [\cite{DBLP:journals/corr/abs-2206-03130, DBLP:journals/tmlr/RuhkopfMDTHL23}]. However, these works often require a complex learning pipeline and do not specifically target budget-limited scenarios.}

{In this work, we seek to address the budget-limited algorithm selection problem from a different perspective. We frame this problem as a special type of Markov Decision Process, paving the way for Reinforcement Learning methods. In this setting, an agent actively requests to train and test algorithms to reveal their performances on a given dataset, which implies an ``active meta-learning'' setting. The agent should interrupt the training of less promising algorithms based on partial learning curves, rather than waiting until all algorithm candidates are fully trained to evaluate. This idea is inspired by the concept of ``freezing'' and ``thawing'' algorithms proposed by \cite{DBLP:journals/corr/SwerskySA14}. Our framework allows meta-learning to improve the algorithm selection process. We focus on First-level meta-learning, which involves algorithm evaluations (i.e., learning curves), dataset meta-features, and/or algorithm hyperparameters; in the same line of work by [\cite{DBLP:conf/pkdd/NguyenGSG21, DBLP:phd/hal/SunHosoya19, DBLP:conf/pkdd/Sun-HosoyaGS18}]. The meta-trained agent should balance two aspects: (1) \textit{exploration-exploitation} trade-offs between continuing to train an already tried good candidate and trying a new candidate; and (2) \textit{multi-fidelity} trade-offs between querying high-fidelity data with higher cost and low-fidelity data with lower cost.} 

We summarize our contributions in this work as follows: 
\begin{enumerate}[(i)]
    \item We formulate the limited-budget algorithm selection problem using Markov Decision Processes (MDP). Our framework neither requires explicit extrapolation of learning curves nor does it necessitate parametric and concavity assumptions. (Section 2)

    \item We discuss our challenge series design, including the novel benchmark datasets we created and used in the challenges. We present the challenge results compared to several baseline methods instantiated from our framework. (Section 3)

    \item We perform a comprehensive result analysis with a comparison of data usage and policy types of methods (Section 3.3.1). To see the benefits of meta-learning and learning from learning curves, we conduct an ablation study for a specific baseline method (Section 4.2). We examine in-depth the strategies learned by the winning methods, comparing them with the baselines and offering diverse effective approaches. (Section 4.3)
\end{enumerate}

\rv{Code for reproducing the challenge results, as well as the analysis results presented in this paper, is publicly available in our repositories.
\footnote{https://github.com/LishengSun/metaLC-post-challenge-analysis-1st-round}
\footnote{https://github.com/LishengSun/metaLC-post-challenge-analysis-2n-round}.}
\begin{figure}
    \centering    \includegraphics[width=0.48\textwidth]{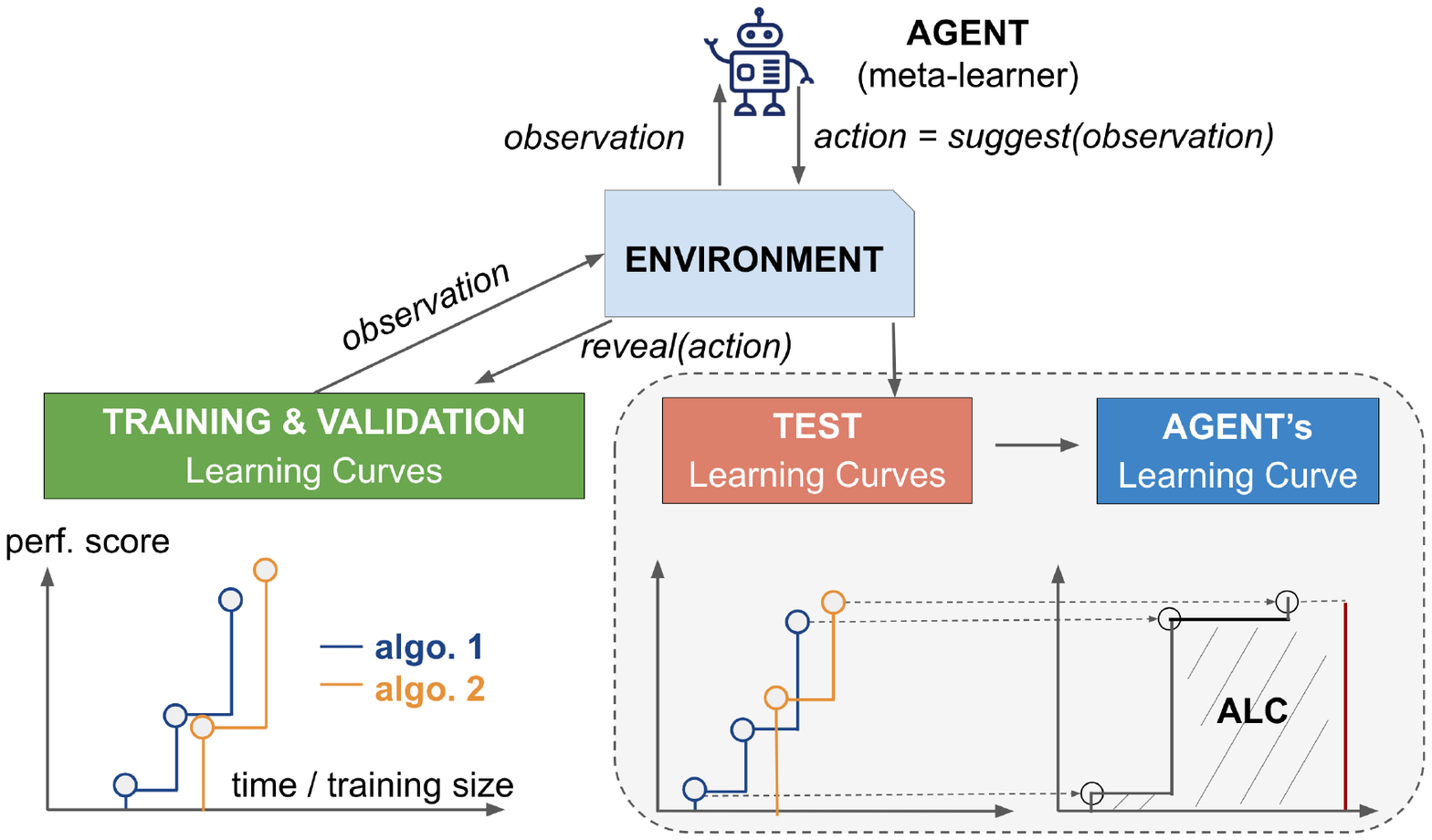}
    \caption{\textbf{Our problem setup}. Given a dataset, an agent (meta-learner) $\mathcal{M}$ takes an action to start or continue training an algorithm using a budget, based on an observation containing partially revealed training and validation learning curves. The corresponding test learning curves are kept hidden and used for computing a reward to be returned to the agent. This interaction is repeated until the given total budget is exhausted. }
    \label{fig:meta-testing}
\end{figure}
\section{Problem Setup}
\label{sec:problem-setup}
{We formally describe the budget-limited algorithm selection problem. Let $\mathcal{D}$ represent a set of datasets, and $\Omega$ denote a set of algorithms of interest. Given a dataset $d_i \in \mathcal{D}$, the goal is to identify the best-performing algorithm $\omega_{j^*} \in \Omega$ for $d_i$. During the selection process, a computational resource budget $\mathcal{T}_i$ is provided for training and evaluating algorithms on dataset $d_i$. We are interested in learning a selection policy $\pi$ that efficiently utilizes budget $\mathcal{T}_i$ to find the best-performing algorithm for any given dataset in $\mathcal{D}$.} 

{We formulate this problem as a Markov Decision Process (MDP) and define its components as follows. In an episode, a dataset $d_i$ is given. At time step $t$, a state ${s_t}$ contains current information about learning curves of algorithms in $\Omega$, on the training set $d_i^{train}$ and validation set $d_i^{val}$, denoted by $\mathcal{L}^{train}_{ij}(\tau_j)$ and $\mathcal{L}^{val}_{ij}(\tau_j)$, respectively. Here, a learning curve is a function of cost $\tau_j$ ($0 \leq \tau_j \leq \mathcal{T}_i$), which returns the performance of algorithm $\omega_j$ on dataset $d_i$ after spending $\tau_j$ for training and evaluation, w.r.t. a certain metric (e.g., classification accuracy). The cost $\tau_j$ can be time spent, the percentage of training data used, etc. Thus, the learning curve of algorithm $\omega_j$ is only partially revealed up until $\tau_j$. The state $s_t$ also contains meta-features of dataset $d_i$ and hyperparameters of algorithms in $\Omega$. An agent $\mathcal{M}$ observes state $s_t$ and takes an action $a_t = (\omega_j, \Delta_j, \widehat{\omega_{j^*}})$. In this tuple, $\omega_j$ is an algorithm to be continually trained and evaluated using a budget increment $\Delta_j$ (hence, we update $\tau_j=\tau_j+\Delta_j$ ), and $\widehat{\omega_{j^*}}$ is the predicted best performing algorithm on the test set $d_i^{test}$ given observations. Once action $a_t$ is executed, learning curves $\mathcal{L}^{train}_{ij}(\tau_j)$ and $\mathcal{L}^{val}_{ij}(\tau_j)$ are updated and revealed to the agent in the next time step. In contrast, learning curves on the test set $\mathcal{L}^{test}_{ij}(\tau_{j^*})$ are kept hidden from the agent and used for computing a reward $r_t$. Concretely, $r_t$ is defined as the improvement on the test set, weighted by the normalized remaining budget: 
\begin{equation}
    r_t = \left[\mathcal{L}^{test}_{ij^*_{t}}(\tau_{j^*_t}) - \mathcal{L}^{test}_{ij^*_{t-1}}(\tau_{j^*_{t-1}})\right] \left[1-\Tilde{t}\right]
    \label{eq:reward}
\end{equation}
with the normalized time:
\begin{equation}
    \Tilde{t} = \frac{\log(1+\sum^{|\Omega|}_{j=1}\tau_j/\sigma)}{\log(1+\mathcal{T}_i/\sigma)}
    \label{eq:reward2}
\end{equation}
where $\sum^{|\Omega|}_{j=1}\tau_j$ denotes the total budget spent up to time step $t$. The hyperparameter $\sigma$ controls the emphasis on performance importance at the beginning of the episode. The goal is to encourage the agent $\mathcal{M}$ to discover good algorithms as quickly as possible. In this way, even if the agent is stopped early, we will get as good performance as possible. This is known as an ``Any-time Learning'' setting. An episode is terminated when the total budget $\mathcal{T}_i$ is exhausted. \rv{By integrating the learning curve using horizontal rectangles, in the style of Lebesgue integrals, the accumulated reward during the episode is equal to the area under the learning curve (ALC) of the agent.} Figure \ref{fig:meta-testing} illustrates our problem setup.}

\rv{The proposed MDP framework was tested as follows. During the meta-training phase, participants could employ any learning method using all available resources, including training, validation, and test learning curves, as well as meta-features for all meta-training datasets, to develop the agent’s policy. This policy determines the probability of taking specific actions under certain circumstances. In the meta-testing phase, the agent, using the developed policy, interacts with the environment. The environment provides a new dataset from the meta-test datasets along with observations containing partially revealed training and validation learning curves based on the agent's actions. The corresponding test learning curves remain hidden within the environment and are used solely to compute the reward returned to the agent, contributing to the calculation of the agent's final ALC score.}

{In this MDP, the agent's action does not influence how underlying states and rewards are generated, which distinguishes it from a standard MDP. Consequently, both states and rewards can be pre-computed at the beginning of an episode. Thus, the problem becomes a REVEAL game, a special MDP discussed in prior works by \cite{DBLP:phd/hal/SunHosoya19} and \cite{DBLP:conf/pkdd/NguyenGSG21}.}

\setlength{\tabcolsep}{6pt}
\begin{table*}
\centering
\caption{\textbf{Data usage and policy type of the top-3 teams in each round vs. five baselines.} The {\color{lightblue}{\cmark}} symbol corresponds to the method that meta-learned from learning curves in meta-training, represented by a {\color{lightblue}{blue bar}} in Figure \ref{fig:challenge-result}. A \textit{combined} policy refers to a combination of learned and hard-coded policies.}
\label{table:comparison}
\scalebox{0.75}{
\begin{tabular}{c||ccc|ccc||cc} 
\toprule
\multicolumn{1}{l||}{}                                         & \multicolumn{6}{c||}{\textbf{Data usage}}                                                                                                                                                                                                                                                                                                                                                                                  & \multicolumn{2}{c}{\textbf{Policy type}}                                                                                                                                              \\ 
\cline{2-9}
\multicolumn{1}{l||}{}                                         & \multicolumn{3}{c|}{in \textit{meta-training }}                                                                                                                                                             & \multicolumn{3}{c||}{in \textit{meta-testing}}                                                                                                                                                               & \multicolumn{1}{l}{\multirow{2}{*}{\begin{tabular}[c]{@{}l@{}}\\algorithm\\selection\end{tabular}}} & \multirow{2}{*}{\begin{tabular}[c]{@{}c@{}}\\budget \\allocation\end{tabular}}  \\ 
\cline{2-7}
                                                               & \begin{tabular}[c]{@{}c@{}}learning curve\\progression\end{tabular} & \begin{tabular}[c]{@{}c@{}}dataset\\meta\\features\end{tabular} & \begin{tabular}[c]{@{}c@{}}algorithm \\hyperparameters\end{tabular} & \begin{tabular}[c]{@{}c@{}}learning curve\\progression\end{tabular} & \begin{tabular}[c]{@{}c@{}}dataset \\meta\\features\end{tabular} & \begin{tabular}[c]{@{}c@{}}algorithm \\hyperparameters\end{tabular} & \multicolumn{1}{l}{}                                                                                &                                                                                 \\ 
\hline
Team MoRiHa                                                    & \textcolor{lightblue}{\textbf{\cmark}}                                                     & \cmark                                                          &                                                                     & \cmark                                                              & \cmark                                                           &                                                                     & combined                                                                                            & combined                                                                        \\
Team neptune                                                   & \textcolor{lightblue}{\textbf{\cmark}}                                                     & \cmark                                                          & \cmark                                                              &                                                                     & \cmark                                                           & \cmark                                                              & combined                                                                                            & combined                                                                        \\
Team AIpert                                                    & \textcolor{lightblue}{\textbf{\cmark}}                                                     & \cmark                                                          &                                                                     &                                                                     & \cmark                                                           &                                                                     & combined                                                                                            & combined                                                                        \\ 
\hline
Team dragon\_bra                                               & \textcolor{lightblue}{\textbf{\cmark}}                                                     & \cmark                                                          &                                                                     &                                                                     & \cmark                                                           &                                                                     & combined                                                                                            & hard-coded                                                                      \\
Team diaprofesser                                              & \textcolor{lightblue}{\textbf{\cmark}}                                                     & \cmark                                                          &                                                                     &                                                                     & \cmark                                                           &                                                                     & combined                                                                                            & combined                                                                        \\
Team carml                                                     & \textcolor{lightblue}{\textbf{\cmark}}                                                     & \cmark                                                          & \cmark                                                              & \cmark                                                              & \cmark                                                           & \cmark                                                              & combined                                                                                            & hard-coded                                                                      \\ 
\hline
\textsc{DDQN} {[}\cite{DBLP:conf/aaai/HasseltGS16}]            & \textcolor{lightblue}{\textbf{\cmark}}                                                     &                                                                 &                                                                     & \cmark                                                              &                                                                  &                                                                     & learned                                                                                             & hard-coded                                                                      \\
\textsc{AvgRank} {[}\cite{DBLP:conf/pkdd/NguyenGSG21}]         &                                                                     &                                                                 &                                                                     &                                                                     &                                                                  &                                                                     & combined                                                                                            & hard-coded                                                                      \\
\textsc{Freeze-Thaw} {[}\cite{DBLP:journals/corr/SwerskySA14}] &                                                                     &                                                                 &                                                                     & \cmark                                                              &                                                                  & \cmark                                                              & hard-coded                                                                                          & hard-coded                                                                      \\
\textsc{BoS} {[}\cite{DBLP:conf/pkdd/NguyenGSG21}]             &                                                                     &                                                                 &                                                                     & \cmark                                                              &                                                                  &                                                                     & hard-coded                                                                                          & hard-coded                                                                      \\
\textsc{RandSearch}                                            &                                                                     &                                                                 &                                                                     &                                                                     &                                                                  &                                                                     & hard-coded                                                                                          & hard-coded                                                                      \\
\bottomrule
\end{tabular}}
\end{table*}

\section{MetaLC Challenge}
{We investigated the potentials of meta-learning from learning curves for improving budget-limited algorithm selection, by creating a series of challenges, namely, \textit{MetaLC}. The main objective of these challenges was to train an agent $\mathcal{M}$ (also referred to as a meta-learner) that is capable of meta-learning from learning curves on other datasets and efficiently identify the best-performing algorithm for a new dataset within a limited budget. The challenge series comprised two rounds. 

\subsection{Benchmark meta-datasets}
\label{sec:data}
At the time of organizing this challenge, there were limited meta-datasets of learning curves available in the machine learning (ML) community, despite their widespread use. To facilitate benchmarking, we created a new meta-dataset comprising learning curves of Automated Machine Learning (AutoML) algorithms on 30 cross-domain AutoML datasets provided by \cite{DBLP:books/sp/19/GuyonSBEELJRSSSTV19}. The application domains of these datasets include medical diagnosis, text classification, customer satisfaction prediction, speech recognition, object recognition. The datasets have been preprocessed into suitable fixed-length vectorial representations. 

{We created a set of algorithms by modifying an AutoML baseline provided in the AutoML challenge by \cite{DBLP:books/sp/19/GuyonSBEELJRSSSTV19}, and varying its hyperparameters.\footnote{https://github.com/ch-imad/AutoMl\_Challenge/blob/master/Starting\_kit} Concretely, we changed only the core algorithm of the AutoML method and kept the rest of its components unchanged. In round 1, we used tree-based algorithms (Random Forest, Gradient Boosting) as the core algorithm. To compute a learning curve, we incrementally increased the number of estimators and repeated evaluation. In round 2, we used Nearest Neighbors, Multilayer Perceptron, Adaboost, and Stochastic Gradient Descent algorithms; we computed a learning curve by increasing the training data size (i.e., from 10\%, 20\%, ..., 100\%). In the first round, we had a total of 600 learning curves (20 algorithms $\times$ 30 datasets). In the second round, we had 1200 learning curves (40 algorithms $\times$ 30 datasets).} In both rounds, meta-features of datasets and hyperparameters of algorithms were also provided. A starter kit with a synthetic meta-dataset created by \cite{DBLP:conf/ijcnn/NguyenSGG22, NguyenCFOL}, was given to participants for practice. 

\subsection{Evaluation protocol}

Participants were asked to develop an agent $\mathcal{M}$ and submit it to be executed on our Codalab competition websites.\footnote{1st round: https://codalab.lisn.upsaclay.fr/competitions/753}\footnote{2nd round: https://codalab.lisn.upsaclay.fr/competitions/4894}. {Each round had two phases: {meta-training} and {meta-testing}. We split $\mathcal{D}$ into $\mathcal{D}^{meta-train}$ and $\mathcal{D}^{meta-test}$, and use learning curves of algorithms in these sets for meta-training and meta-testing agent $\mathcal{M}$, respectively. In {meta-training}, participants were allowed to use any kind of learning on provided data.} In {meta-testing}, the agent $\mathcal{M}$ was evaluated by interacting with an environment as described in Section \ref{sec:problem-setup} and illustrated in Figure \ref{fig:meta-testing}. {The ranking of the agent on the leaderboard was determined by the average accumulated reward over datasets in $\mathcal{D}^{meta-test}$.} 

\subsection{Baseline methods}

\rv{In this section, we introduce our five baseline methods. These are methods we instantiated following our proposed framework and leveraging techniques from existing work to serve as baselines. We then compare these baselines with the methods submitted by participants in the challenges to assess whether any improvements over the baselines were achieved. Each baseline method represents a distinct solution within the proposed framework, developing its own policies for algorithm selection and budget allocation.}

\textbf{\textsc{Double Deep Q-Network (DDQN)}}. As the challenge is inspired by RL, we wanted to evaluate the meta-learning capabilities of RL methods. We started with DDQN, a classic RL method proposed by \cite{DBLP:conf/aaai/HasseltGS16}. We used the data given in the meta-training phase to create an RL environment for training (same setup discussed in Section \ref{sec:problem-setup}). {The agent learned a policy $\pi_{\theta,\theta'}$ with two networks, one for action selection parameterized by $\theta$, the other for value estimation parameterized by $\theta'$. The parameters were updated by minimizing the following loss using trajectories sampled from a replay buffer $\mathbb{B}$: $\theta, \theta' = \argmin_{\theta,\theta'}\mathbb{E}_{(s_t, a_t, r_t, s_{t+1})\sim \mathbb{B}}[(y_t-Q(s_t,a_t;\theta))^2]$ with learning target $y_t$ defined by: $y_t = r_t + \gamma Q(s_{t+1}, \argmax_{a}Q(s_{t+1}, a;\theta);\theta')$. We used this learned policy only for choosing algorithm $\omega_j$ in the action triplet $(\omega_j, \Delta_j, \widehat{\omega_{j^*}})$, i.e. $\pi_{\theta,\theta'}(s_t) = \omega_j$. The algorithm with the highest performance revealed on the validation set so far was selected as the predicted best-performing algorithm $\widehat{\omega_{j^*}}$ with $j^* = \argmax_{j}\mathcal{L}^{val}_{ij}(\tau_j)$. For budget allocation, it used a fixed policy that starts with a pre-defined small budget and doubles the amount of budget spent for an algorithm $\omega_j$ every time it resumes training: $\tau_j = 2*\tau_j$.}

\textbf{\textsc{Freeze-Thaw}}. {We considered the problem as a hyperparameter search, as our dataset of algorithms was created by varying their hyperparameters (see Section \ref{sec:data}). We used the Freeze-Thaw Bayesian Optimization (\textsc{Freeze-Thaw}) method proposed by \cite{DBLP:journals/corr/SwerskySA14}, which was adapted to become one of the winning solutions in the AutoML challenge [\cite{DBLP:books/sp/19/GuyonSBEELJRSSSTV19}]. It used the partial learning curve information in a Bayesian Optimization fashion to decide which algorithm to ``freeze'' or ``thaw'' training it. Concretely, it selected algorithm $\omega_j$ that maximizes an acquisition function, i.e. $j = \argmax_{j}f_{acq}(j)$. It used an entropy search acquisition function that maximizes the expected information gain over the location of the performance maximum: $f_{acq}(j)= \int(H(P^y_{max}) - H(P_{max}))P(y|\{(j^n,y^n)\}_{n=1}^N)dy$. Here, $P_{max}$ represents the current estimated distribution over the performance maximum, and $P^y_{max}$ is the updated distribution given that point $j$ yields the performance $y$. $N$ is the number of observations. Similar to \textsc{DDQN}, $\widehat{\omega_{j^*}}$ was set by $j^*= \argmax_j\mathcal{L}^{val}_{ij}(\tau_j)$. We used a fixed budget increment $\Delta_j$ in each time step. We note that \textsc{Freeze-Thaw} does not have a meta-learning capability.}

\textbf{\textsc{AvgRank}} {Inspired by existing works [\cite{Abdulrahman2018, DBLP:conf/ecml/BrazdilS00, DBLP:conf/mldm/LeiteBV12,Lin2010}], this \textsc{AvgRank} baseline meta-learned an average ranking of algorithms $\Omega$ on $\mathcal{D}^{meta-train}$ during the meta-training phase. In meta-testing, only the algorithm $\omega_j$ that ranked highest was chosen to be trained and evaluated using the entire given budget $\Delta_j = \mathcal{T}_i$. The average rank of algorithm $\omega_j$ was defined as: $avg\_rank({\omega_j})= \left[\sum_{d_i \in \mathcal{D}^{meta-train}} rank(\omega_j, d_i)\right]/|\mathcal{D}^{meta-train}|$. Since this baseline trained only one algorithm, hence, $\widehat{\omega_{j^*}}=\omega_j$. In a real-life scenario, this is an expensive baseline, as it requires training and testing the entire set of algorithms on all meta-training datasets to have an accurate ranking.}

\textbf{\textsc{BestOnSamples (BoS)}}. {This baseline, inspired by the work of \cite{Petrak2000}, selected the algorithm that performed best within a fixed small budget $\alpha$. More specifically, at the beginning of each episode, it trained every algorithm with the same budget $\Delta_j = \alpha$. Based on the observed results, it selected algorithm $\omega_j$ that achieved the highest performance with $j = \argmax_{j} \mathcal{L}_{ij}^{val}(\alpha)$,
and spent the entire remaining budget (now, $\Delta_j = \mathcal{T}_i - \alpha*|\Omega|$). Again, the predicted best-performing algorithm $\widehat{\omega_{j^*}}$ was set by $j^*= \argmax_{j}\mathcal{L}^{val}_{ij}(\tau_j)$. This baseline does not have a meta-learning capability.} 

\textbf{\textsc{RandSearch}}. { This simple baseline performs a random search over the given algorithm set $\Omega$. It uniformly sampled an algorithm $\omega_j \sim \mathcal{U}(\Omega)$ and uniformly assigned an amount of budget for training and evaluating the algorithm, $\Delta_j \sim \mathcal{U}(\Delta_{min}, \Delta_{max})$. Again, $\widehat{\omega_{j^*}}$ is set by $j^* = \argmax_{j}\mathcal{L}^{val}_{ij}(\tau_j)$. Due to its high variance, we ran this baseline five times internally and reported its average performance. This baseline was implemented solely for comparison purposes and is not realistic, as one would not average several runs of an algorithm in practice.}

\begin{figure*}
    \centering
    \begin{subfigure}[b]{0.245\textwidth}
        \centering
        {\includegraphics[width=1\linewidth]{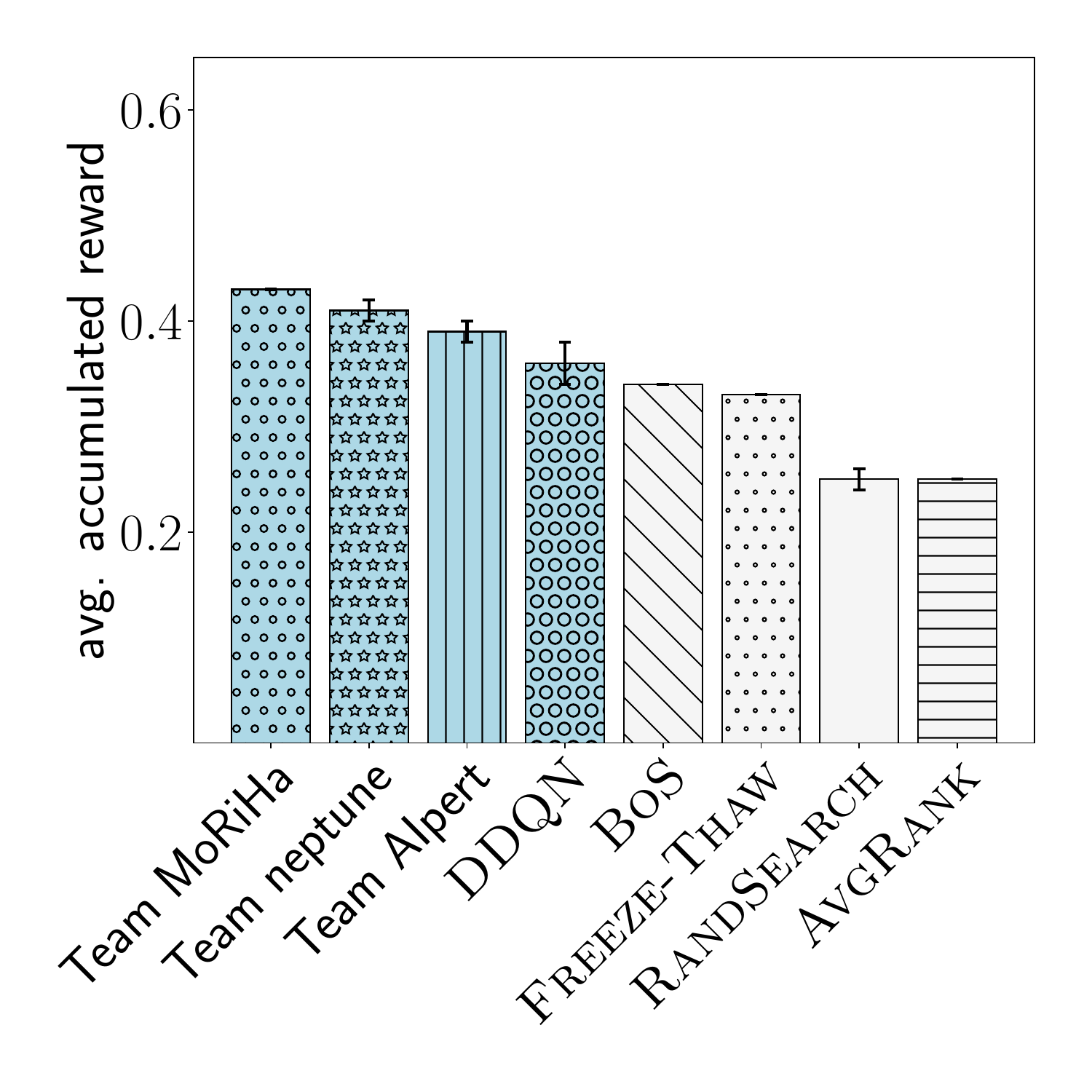}}
        \caption{{1st Round - Any-time Learning}}
    \end{subfigure}
    \hfill
    \begin{subfigure}[b]{0.245\textwidth}
        \centering
        {\includegraphics[width=1\linewidth]{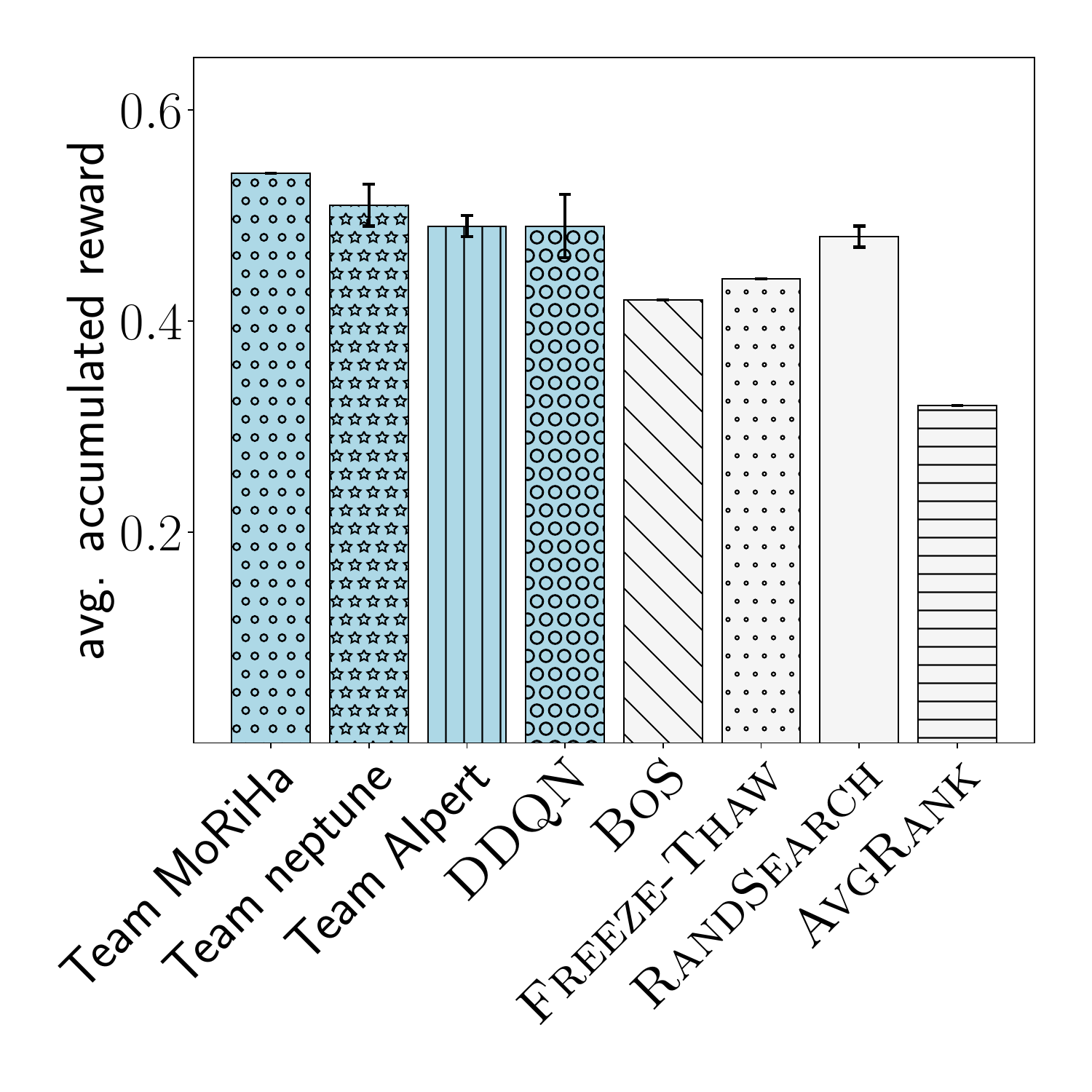}}
        \caption{{1st Round - Fixed-time Learning}}
        \label{fig:1st-roundb}
    \end{subfigure}
    \hfill
    \begin{subfigure}[b]{0.245\textwidth}
        \centering
        {\includegraphics[width=1\linewidth]{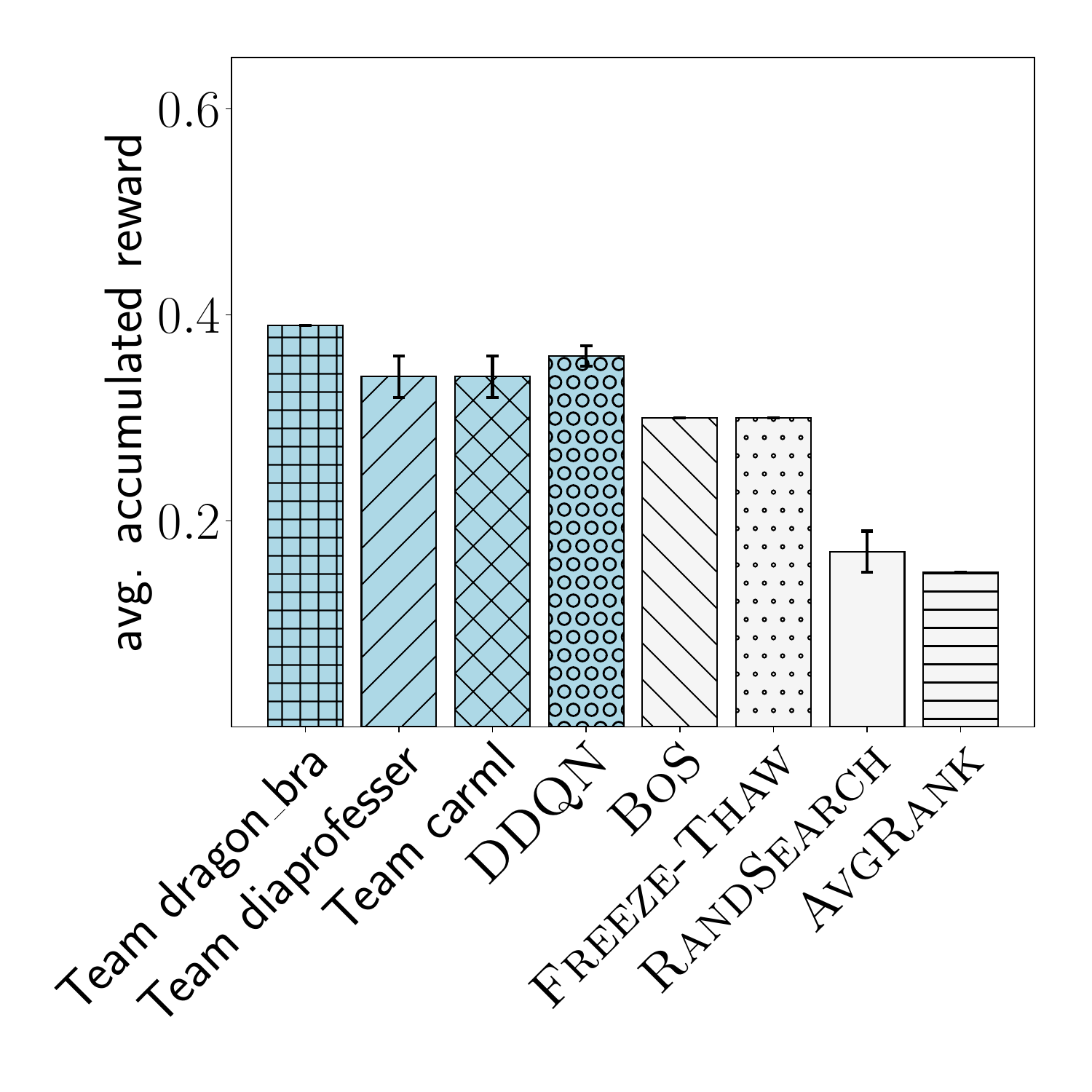}}
        \caption{{2nd Round - Any-time Learning}}
        \label{fig:2nd-rounda}
    \end{subfigure}
    \hfill
    \begin{subfigure}[b]{0.245\textwidth}
        \centering
        {\includegraphics[width=1\linewidth]{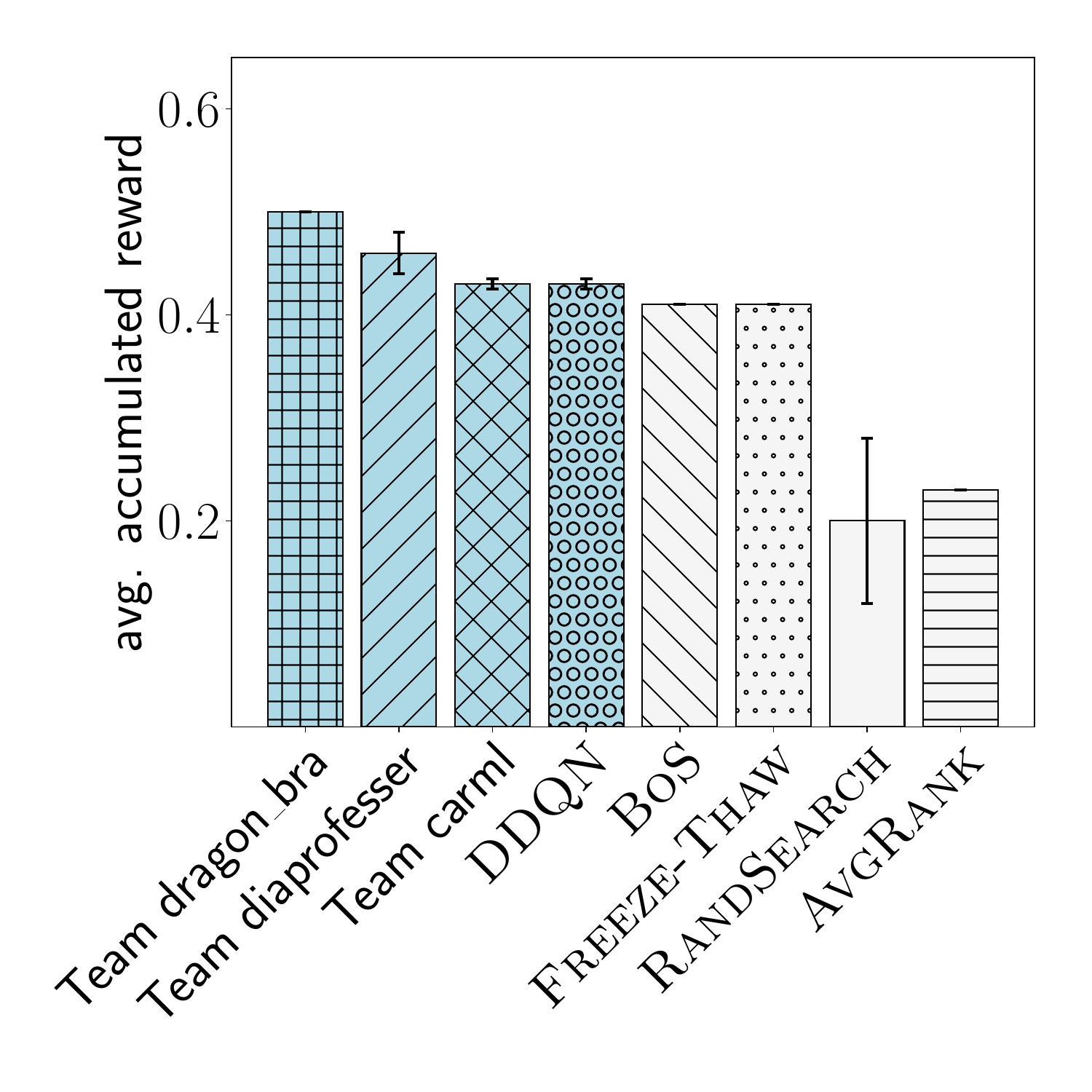}}
        \caption{{2nd Round - Fixed-time Learning}}
        \label{fig:2nd-roundb}
    \end{subfigure}
    \caption{\textbf{MetaLC challenge results.} Comparison of top-3 teams and five baselines. Blue bars represent methods that meta-learned from learning curves in meta-training (corresponds to \textcolor{lightblue}{\textbf{\cmark}} in the first column of Table \ref{table:comparison}). We highlight \textsc{RandSearch} in plain gray, a special baseline with internally averaged performance over several runs. Results for fixed-time learning are included for analysis purposes only and were not officially used in our challenges. The reported results are from the worst run out of three runs with different seeds, and the error bar indicates the standard deviation across meta-test datasets.} 
    \label{fig:challenge-result}
\end{figure*}

\subsubsection{Results}
Table \ref{table:comparison} provides a summary of winning teams' methods and baselines. We found that learning curves and dataset meta-features were more frequently utilized than algorithm hyperparameters. In each round, only one team made use of partially revealed learning curves during meta-testing. Surprisingly, only a few participants employed RL to train their agents, despite the challenge setting being designed as an RL problem. Figure \ref{fig:challenge-result} presents the challenge results, showing the average accumulated reward of each method over meta-test datasets for both the ``Any-time learning'' and ``Fixed-time learning'' settings. In the ``Fixed-time learning'' setting, the accumulated reward is equivalent to the highest performance found during an episode, regardless of when it is found.

In the first round at \textit{WCCI 2022}, the top-3 teams outperformed the best baseline \textsc{DDQN} in the Any-time learning setting; while in the Fixed-time learning setting, \textsc{DDQN} slightly surpassed the third-ranked team. These teams used a combination of learned and hard-coded policies for algorithm selection and budget allocation (Table \ref{table:comparison}). Only team ``MoRiHa'', which ranked first, utilized partially revealed learning curves during meta-testing and achieved the highest ALC score on 21/30 datasets. Notably, they performed better than other teams in some multi-label and multi-class classification tasks, such as \textit{tania}, \textit{robert}, \textit{newsgroups}, and \textit{marco}. One of their key findings was that switching the explored algorithm more than once is rarely beneficial due to the high cost associated.

In the second round at \textit{AutoML-Conf 2022}, among the top-3 teams, only one team (``dragon\_bra'') beat \textsc{DDQN} in Any-time learning, while two teams (``dragon\_bra'' and ``diaprofesser'') outperformed \textsc{DDQN} in Fixed-time learning. All winners meta-learned from learning curves in meta-training. A combined policy for algorithm selection was used by all winners, while 2/3 winners used hard-coded rules for distributing the given budget. Team ``dragon\_bra'' obtained the highest average ALC by winning in 6/15 datasets. They had two key findings: (i) At the beginning of an episode, the agent should spend only a small budget to grab a ``base score'', which would make the area under the agent's learning curve larger; and (ii) if a small budget is used near the end of an episode, it is unlikely that substantial improvement will be made.

\section{Discussions}
Our series of challenges was the first of its kind in the ML community when we introduced it. We delve into important considerations during benchmark creation and challenge design, along with notable observations and investigations.
\subsection{Challenge design}

\textbf{Learning Curve}. Our benchmark meta-datasets contained pre-recorded learning curves to avoid any on-the-fly computational issues during the challenges. In ML community, there are various types of learning curves, such as Iteration Learning Curve, and Observation Learning Curve to name a few; see more types in prior work by \cite{DBLP:journals/corr/abs-2201-12150}. As mentioned in Section \ref{sec:data}, we leveraged the AutoML challenge to collect the learning curves. In the first round, a learning curve was represented as a function of time. The points on the learning curves were irregularly spaced and chosen by the learning algorithms themselves (as set by the AutoML challenge). When a participant's method requested points between two recorded points, we interpolated the requested points using the closest previously recorded point. However, this approach did not provide new information on the learning curve while still incurring a cost. To address this, we introduced a new type of learning curve based on training data size in the second round. Participants' methods could choose a training data size from a fixed given set, query an algorithm's performance w.r.t. the chosen size, and pay the associated cost.

\textbf{Data splitting}. In the first round, we kept validation and test learning curves separate, using them for the feedback phase and final testing phase, respectively. This was to prevent overfitting on the test learning curves. However, it was pointed out in the first round that the test learning curves were highly correlated with the validation curves. Therefore, one could overfit the former by simply overfitting the latter. To mitigate this issue, we divided our dataset into two equal portions in the second round, using one half for the feedback phase and keeping the other half completely private for the final testing phase.

\subsection{{Ablation study of \textsc{DDQN} baseline}}
\begin{figure}
    \centering
    \includegraphics[width=0.41\textwidth]{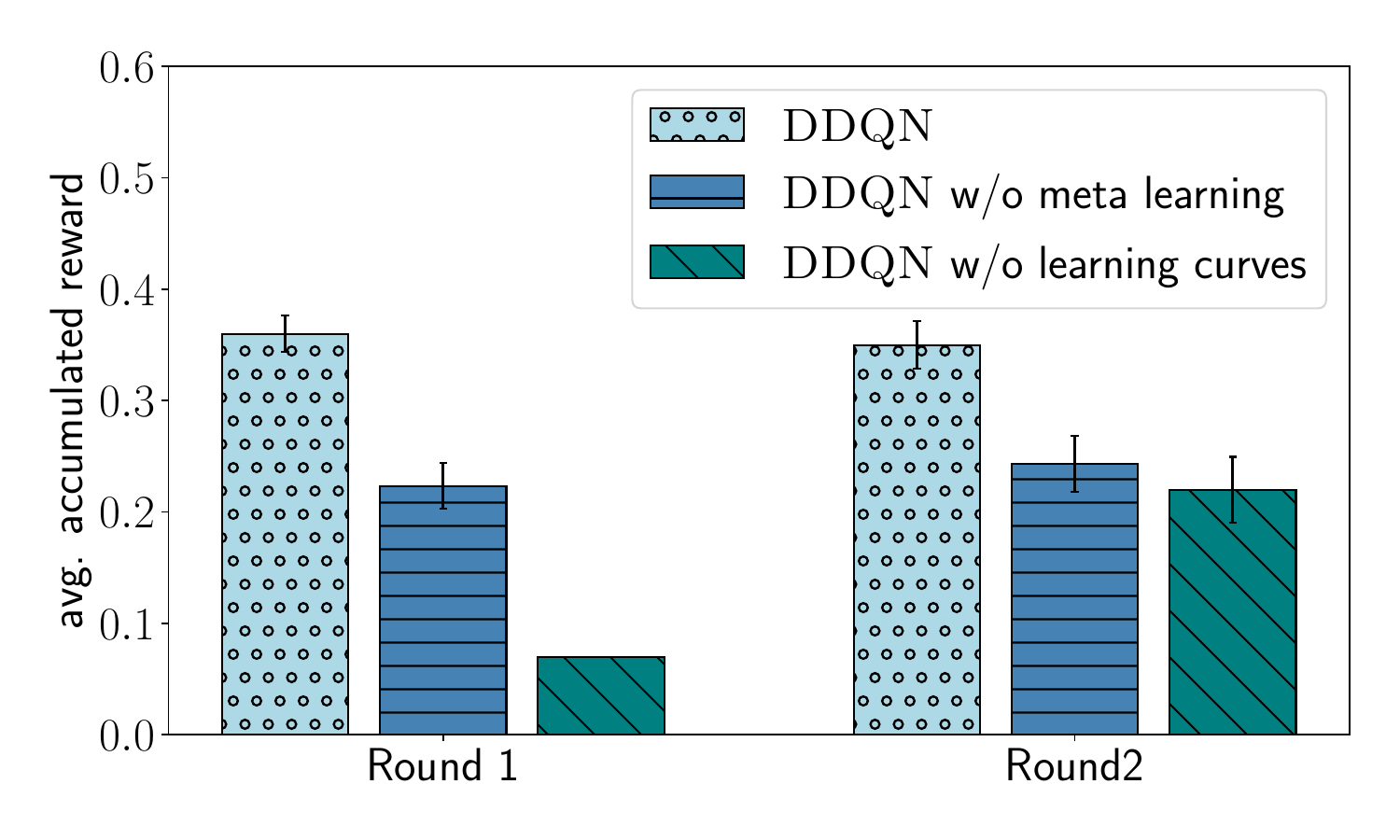}
    \caption{{\textbf{Ablation study of \textsc{DDQN} baseline.} Meta-learning and progression of learning curves improved \textsc{DDQN}'s performance in both challenge rounds.}}
    \label{fig:ablation}
\end{figure}

\begin{figure*}
    \centering
    \begin{subfigure}[b]{0.245\textwidth}
        \centering
        {\includegraphics[width=1\textwidth]{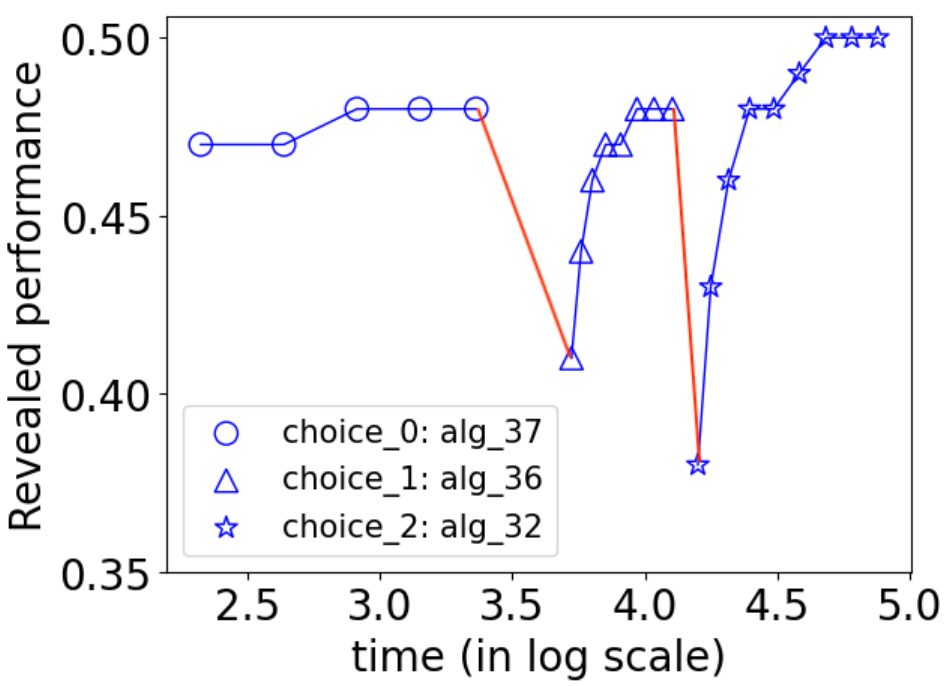}}
        \caption{DDQN}
    \end{subfigure}
    \begin{subfigure}[b]{0.245\textwidth}
        {\includegraphics[width=1\textwidth]{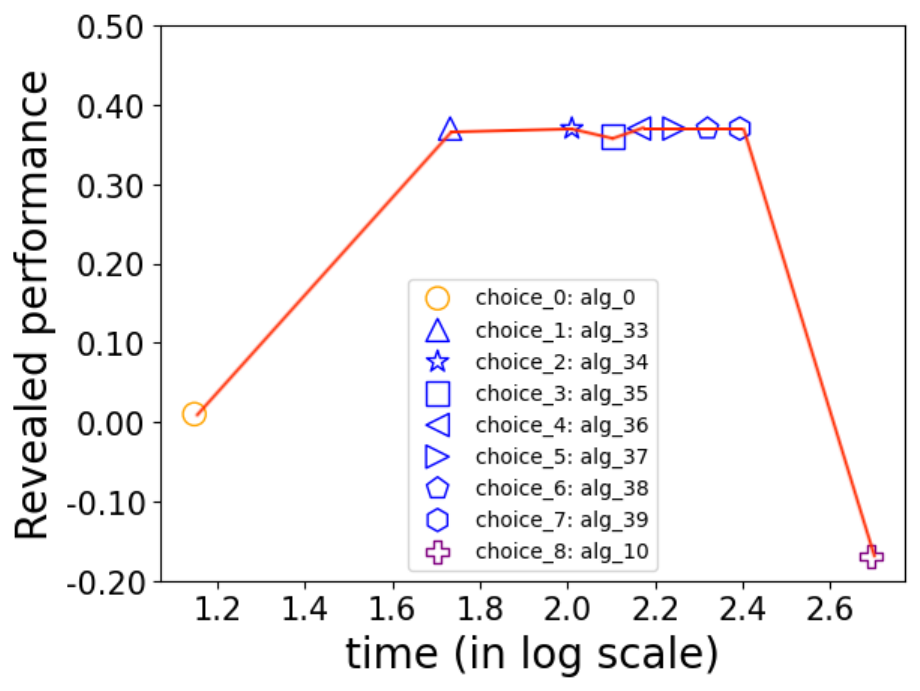}}
        \caption{Team diaprofesser}
    \end{subfigure}
    \begin{subfigure}[b]{0.245\textwidth}
        {\includegraphics[width=1\textwidth]{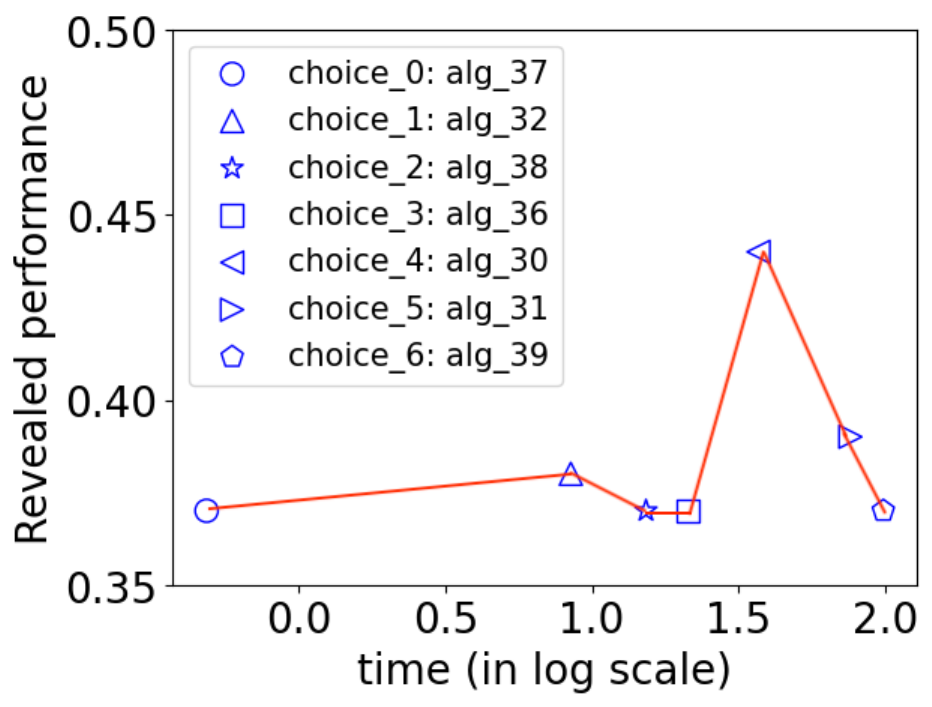}}
        \caption{Team carml}
    \end{subfigure}
    \begin{subfigure}[b]{0.245\textwidth}
        {\includegraphics[width=1\textwidth]{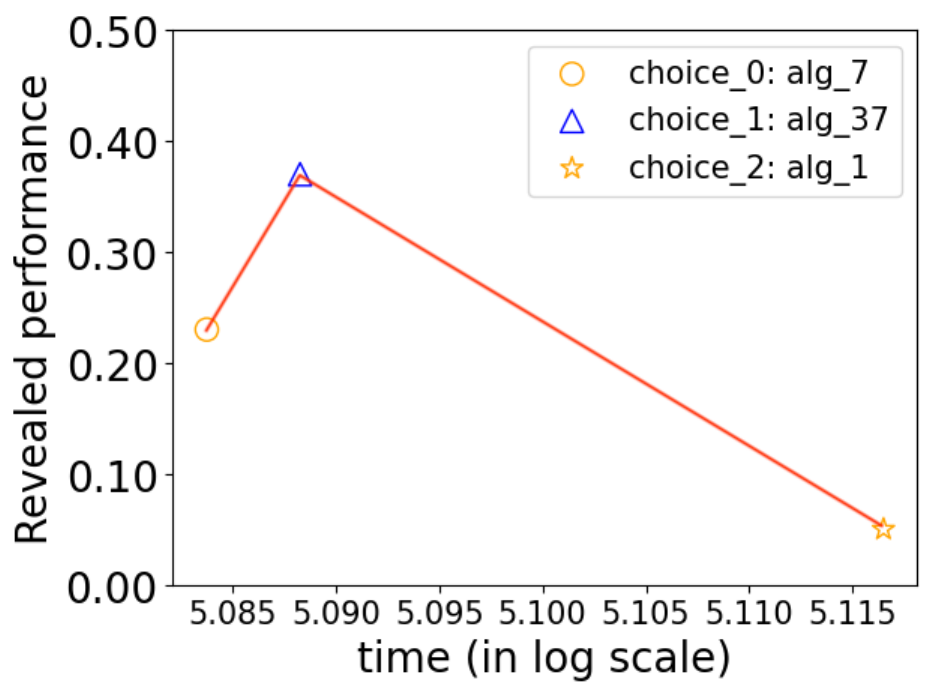}}
        \caption{{Team dragon\_bra}}
    \end{subfigure}
    \caption{{\textbf{Trajectories of baseline DDQN and winning teams' methods on dataset \textit{Flora}, capturing moments of algorithm transitions.} Each marker corresponds to a choice of algorithm $\omega_j$, with the chosen algorithm's family denoted by the marker's color: SGD in blue, AdaBoost in orange, and KNN in purple. Transitions between algorithms are marked with red lines. (a) The DDQN agent began with a strong candidate and consistently selected it. It made a transition only when a performance plateau was reached. (b, c, d) Winning teams' agents exhibited less repetition in their choices and placed a greater emphasis on exploration to achieve better results. The different time ranges on the x-axis were chosen near the beginning of the episode, specifically targeting moments of transition.}}
    \label{fig:ddqn_traj}
\end{figure*}

{We conducted an ablation study to examine the benefits of using learning curves and meta-learning in tackling the budget-limited algorithm selection problem. Figure \ref{fig:ablation} shows a performance comparison of the \textsc{DDQN} baseline with certain components removed. First, we compared the performance of \textsc{DDQN} with and without meta-learning by omitting the meta-training phase. In the latter case, the policy network $\pi_{\theta,\theta'}$ of \textsc{DDQN} was randomly initialized and did not go through any tuning using meta-learning. Second, we compared the performance of \textsc{DDQN} with learning from the entire progression of learning curves (all points on the learning curves) versus learning solely from the final evaluations (last points on the learning curves). In the latter case, an agent selects only $\omega_j$,  and the last points on the pre-computed learning curves $\mathcal{L}_{ij}^{train,val}$ are returned, with an associated cost $\Delta_j$ the agent must pay. The results showcased that both meta-learning and learning curve progression highly contributed to the performance of \textsc{DDQN}. This finding also explains the success of the top-3 teams in both challenge rounds,  as their methods meta-learned from learning curve progression.}

\subsection{Agent policy and action trajectory}
{We investigated the learned policies of agents by examining their action trajectories during the meta-testing phase. Figure \ref{fig:ddqn_traj} shows an example dataset, \textit{Flora}, where we compare trajectories generated by our best baseline, \textsc{DDQN}, and the top-3 teams' methods.  We focus on segments that emphasize moments of algorithm transition, particularly near the beginning of the episode. (1) DDQN's policy: The agent started with algorithm no. 37, a top-performing candidate in meta-training (ranked 3rd/40 algorithms in terms of average performance). It then transitioned to train a new algorithm, no. 36 (ranked 6th), and ultimately to algorithm no. 32 (ranked 7th) upon observing a performance plateau. This pattern demonstrated DDQN's ability to meta-learn effective initial candidates, ensuring a strong starting position. Subsequently, it switched to a new candidate or resumed training a previously paused candidate when the current one reached a performance plateau, guaranteeing efficient any-time learning for which it was meta-trained. (2) Winning teams' policies: In contrast to DDQN, the winning teams' policies prioritized optimal any-time learning performance through a strategy characterized by less repetition of choices and a greater emphasis on exploration. Specifically, team ``Dragon\_bra'' is noted for its cautious approach, thoroughly assessing algorithms, particularly its initial choice, and choosing from different algorithm families like AdaBoost (algorithm no. 7 and no. 1 in orange) to SGD (algorithm no. 37 in blue). Meanwhile, team ``carml'' and team ``Diaprofesser'' favored a more exploratory approach within the most promising regions according to average performance in meta-training (SGD family, algorithm no. 30 to 39, shown in blue markers). This indicates a targeted exploration strategy within these high-potential areas. For a comprehensive view, the complete trajectories for all datasets and a heatmap illustrating algorithm rankings (aiding in identifying promising candidates) can be found in our public repository.}\footnote{https://github.com/LishengSun/metaLC-post-challenge-analysis-2nd-round/tree/main/figs/ddqn\_trajectory}

We observed that a simple method like \textsc{Best on Samples (BoS)} showed decent performance and occasionally surpassed DDQN when algorithms' learning curves do not intersect often, i.e. algorithm ranking does not change significantly w.r.t. the budget spent. Figure \ref{fig:algo-ranking} shows datasets where this phenomenon was observed. More concretely, \textsc{BoS} beat \textsc{DDQN} in 9/30 and 5/15 datasets in the first round and the second round, respectively. This method can serve as a cost-effective yet competitive baseline solution, especially in situations where the implementation of meta-learning is overly complex or resource-intensive.

\begin{figure*}
    \centering
    \begin{subfigure}[b]{0.245\textwidth}
        \centering
        {\includegraphics[width=1\linewidth]{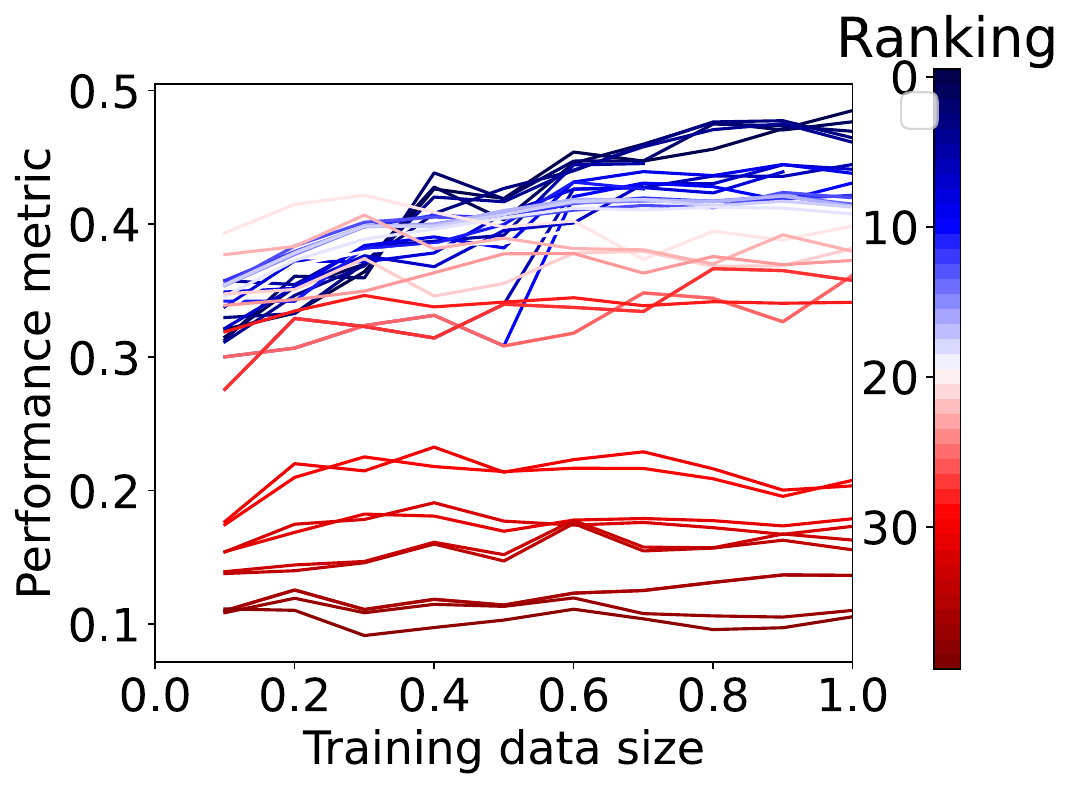}}
        \caption{{waldo}}
        
    \end{subfigure}
    \hfill
    \begin{subfigure}[b]{0.245\textwidth}
        \centering
        {\includegraphics[width=1\linewidth]{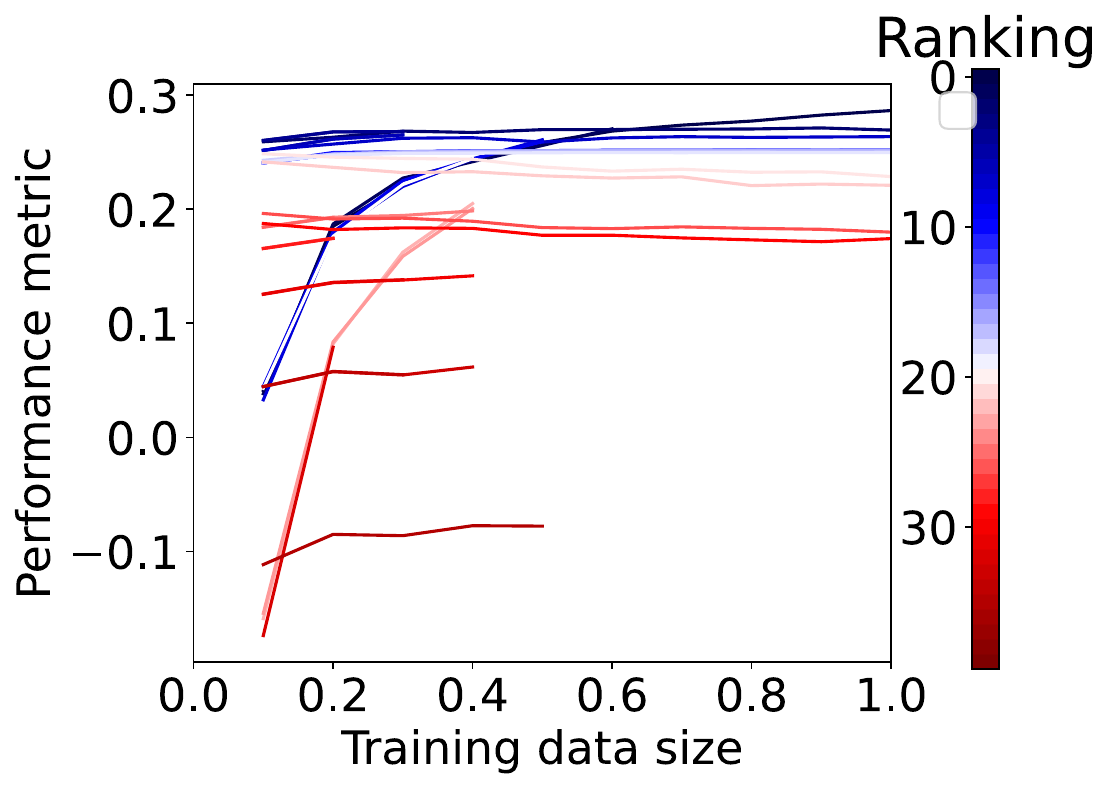}}
        \caption{{pablo}}
        
    \end{subfigure}
    \hfill
    \begin{subfigure}[b]{0.245\textwidth}
        \centering
        {\includegraphics[width=1\linewidth]{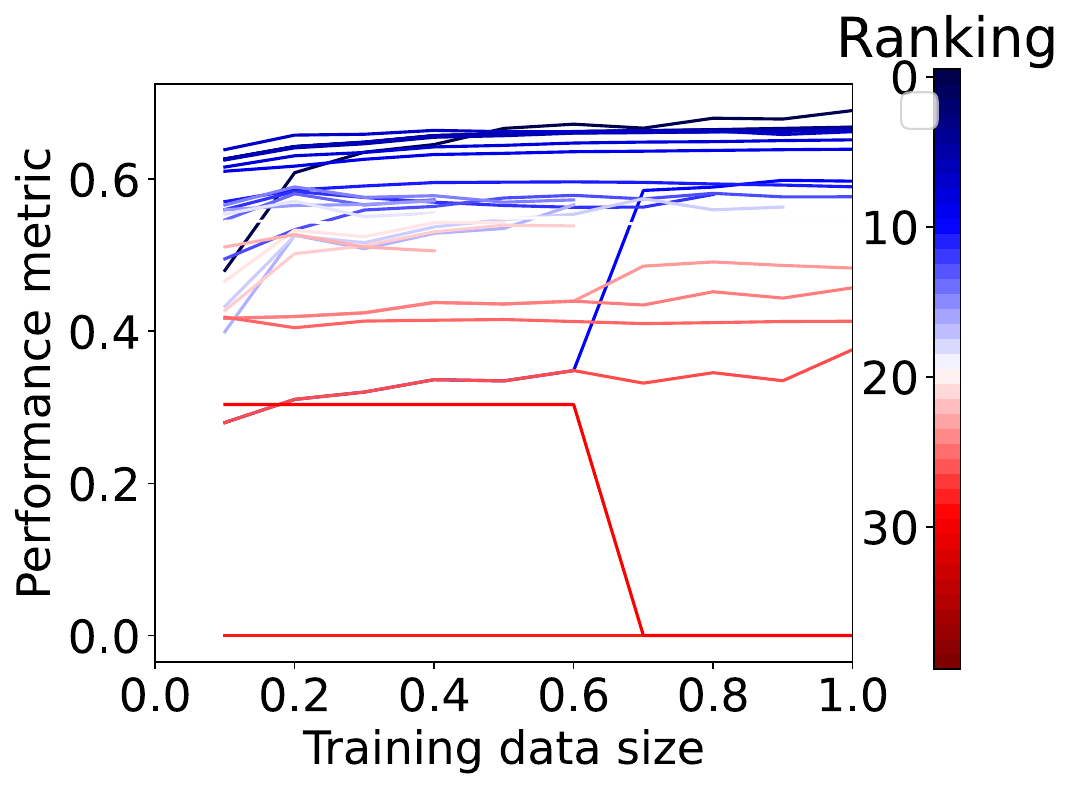}}
        \caption{{marco}}
       
    \end{subfigure}
    \hfill
    \begin{subfigure}[b]{0.245\textwidth}
        \centering
        {\includegraphics[width=1\linewidth]{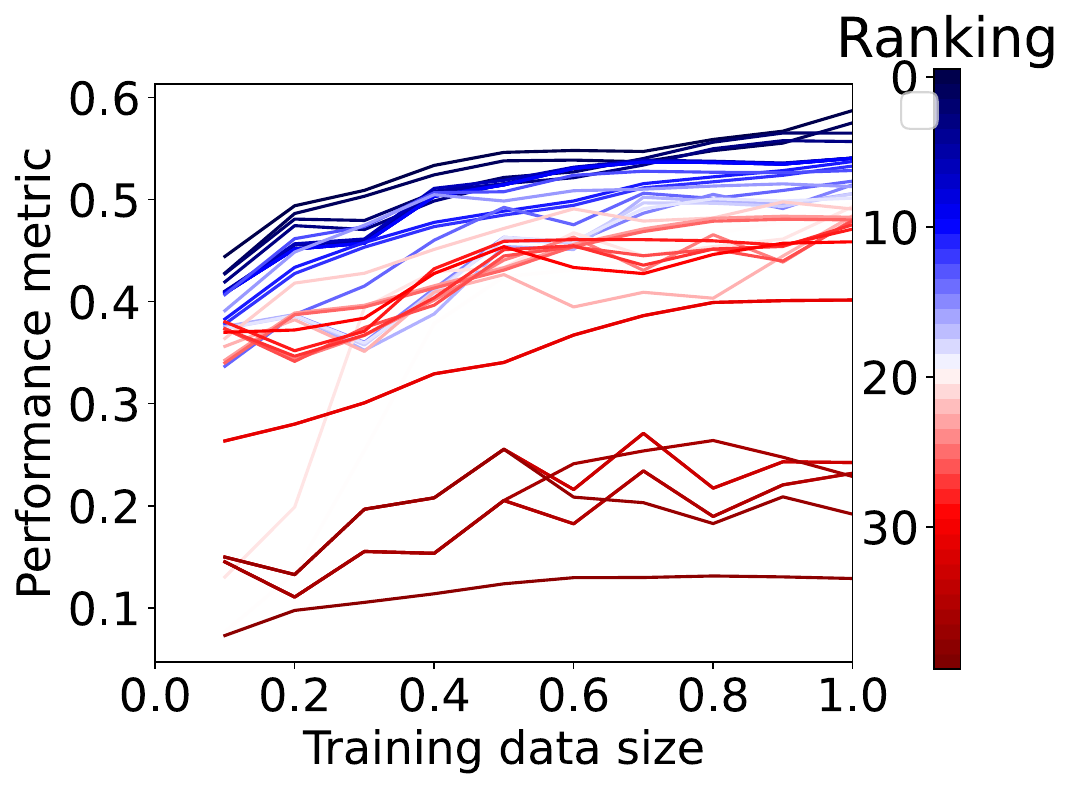}}
        \caption{{evita}}
      
    \end{subfigure}
    
    \begin{subfigure}[b]{0.245\textwidth}
        \centering
        {\includegraphics[width=1\linewidth]{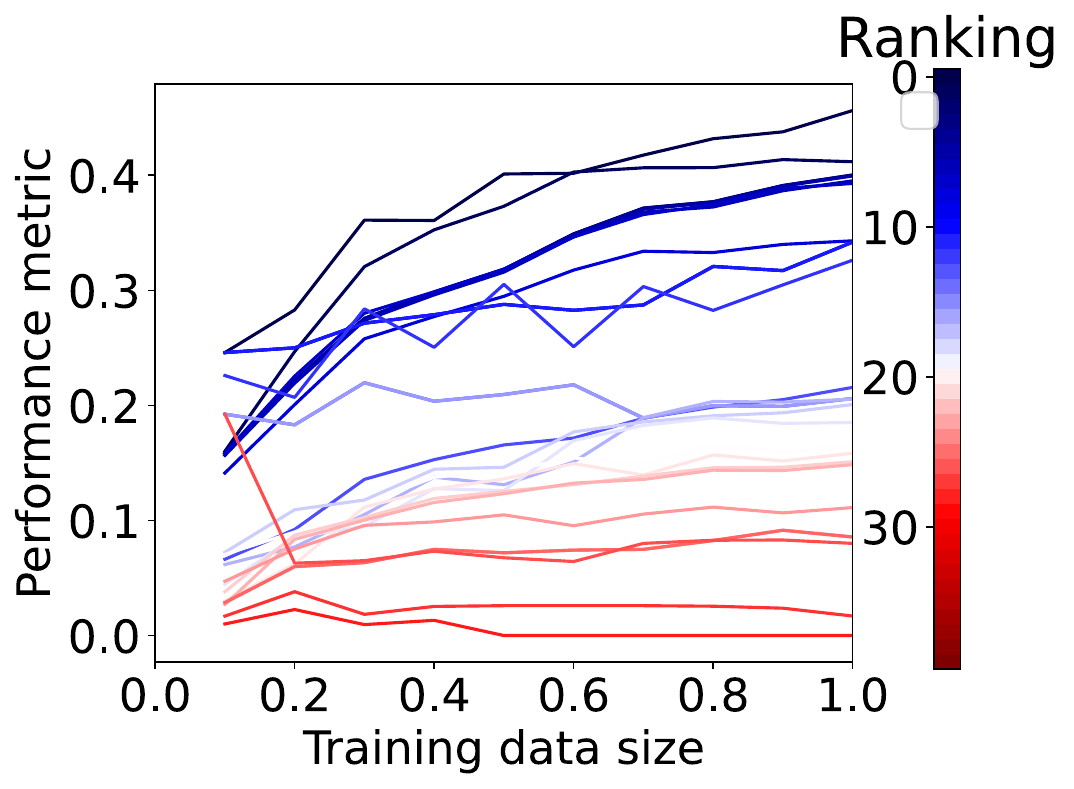}}
        \caption{{wallis}}
        
    \end{subfigure}
    \hfill
    \begin{subfigure}[b]{0.245\textwidth}
        \centering
        {\includegraphics[width=1\linewidth]{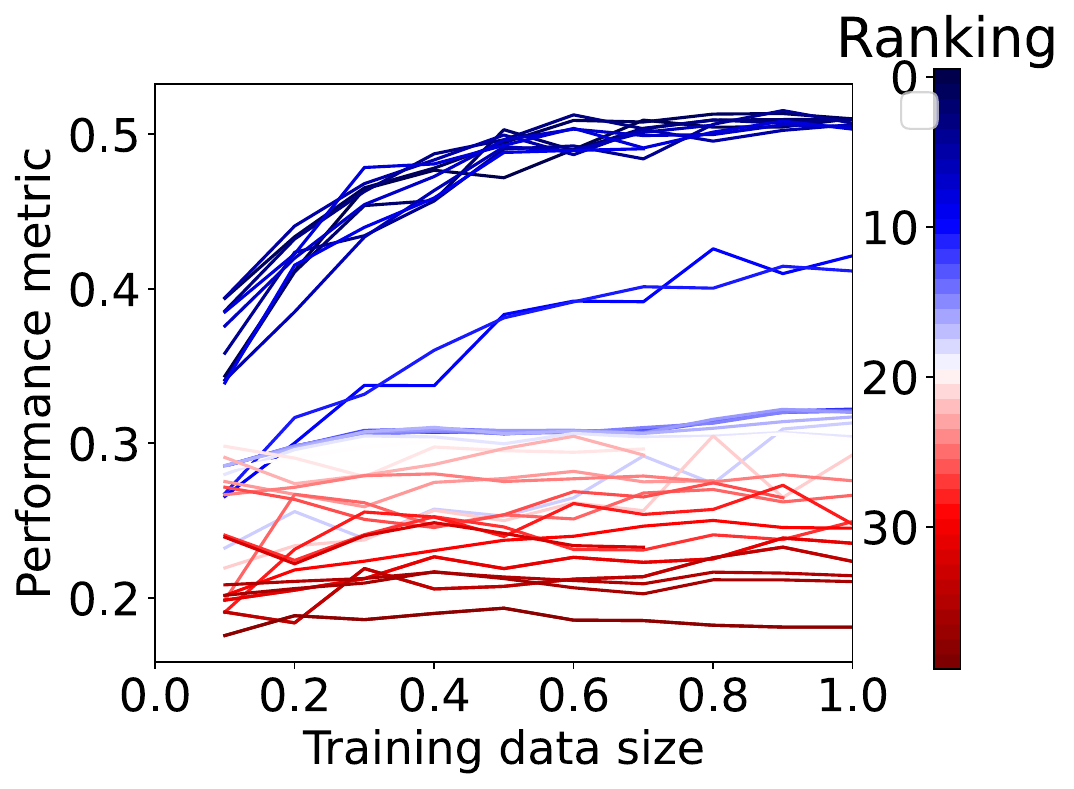}}
        \caption{{jannis}}
       
    \end{subfigure}
    \hfill
    \begin{subfigure}[b]{0.245\textwidth}
        \centering
        {\includegraphics[width=1\linewidth]{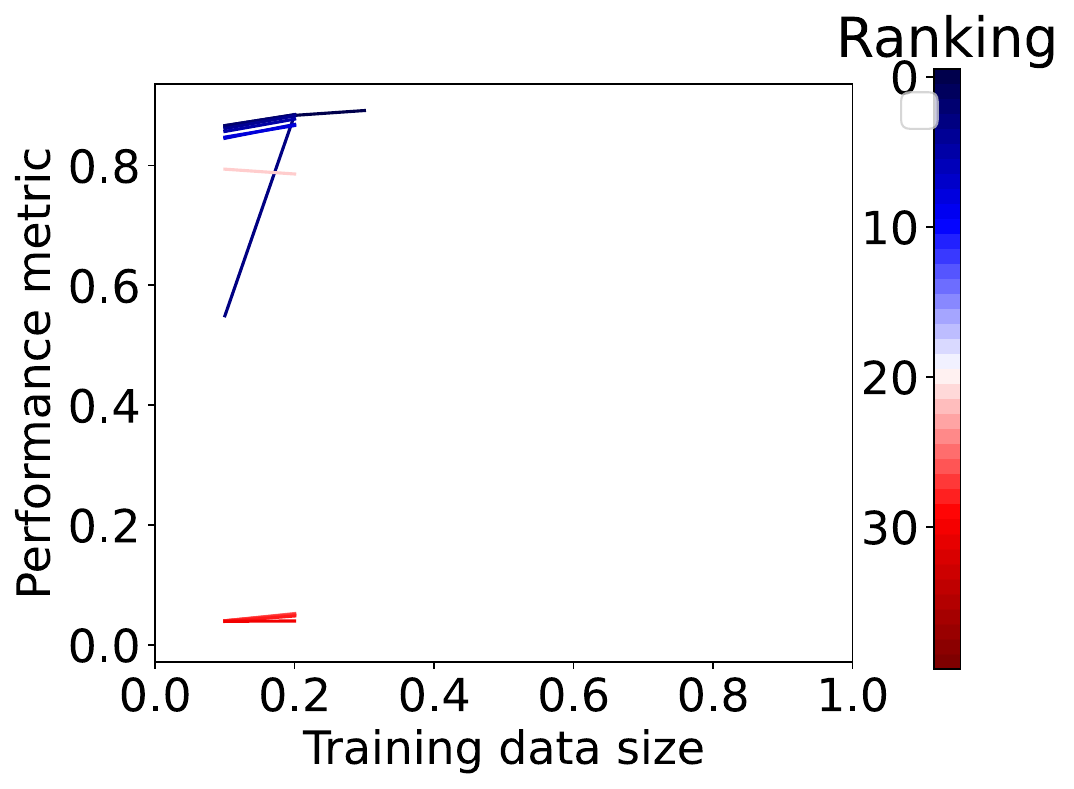}}
        \caption{{dionis}}
       
    \end{subfigure}
    \hfill
    \begin{subfigure}[b]{0.245\textwidth}
        \centering
        {\includegraphics[width=1\linewidth]{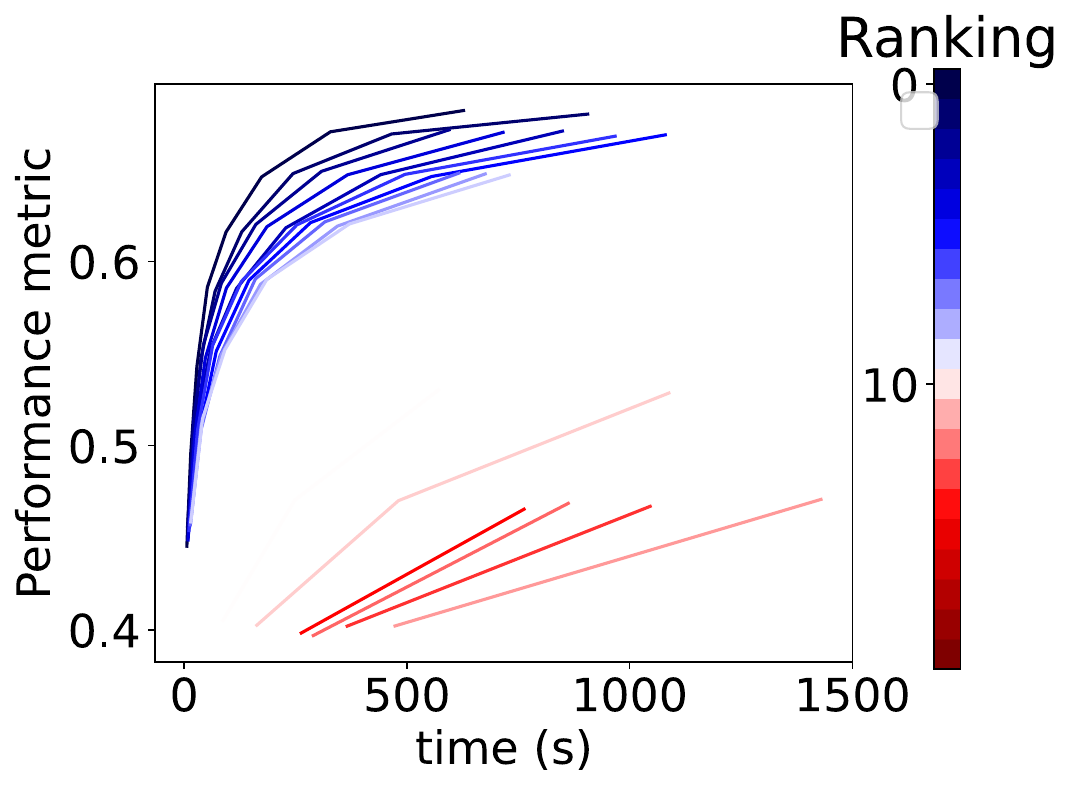}}
        \caption{{alexis}}
       
    \end{subfigure}

    \begin{subfigure}[b]{0.245\textwidth}
        \centering
        {\includegraphics[width=1\linewidth]{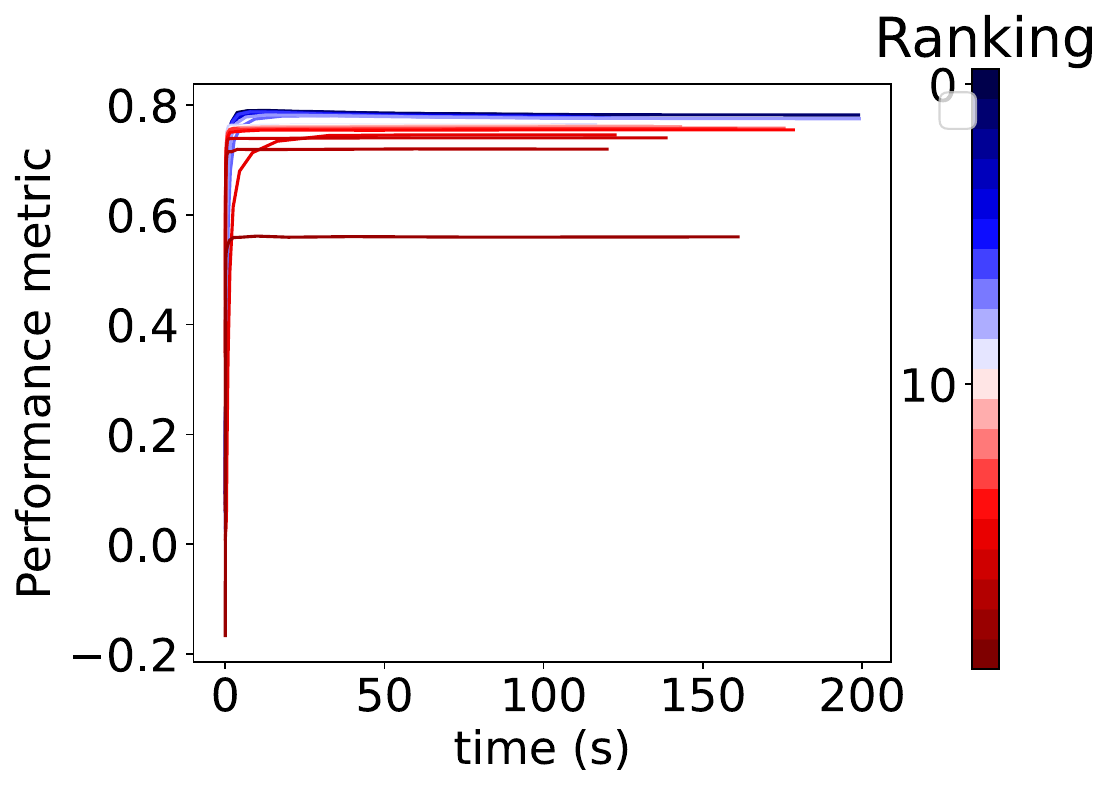}}
        \caption{{cadata}}
       
    \end{subfigure}
    \hfill
    \begin{subfigure}[b]{0.245\textwidth}
        \centering
        {\includegraphics[width=1\linewidth]{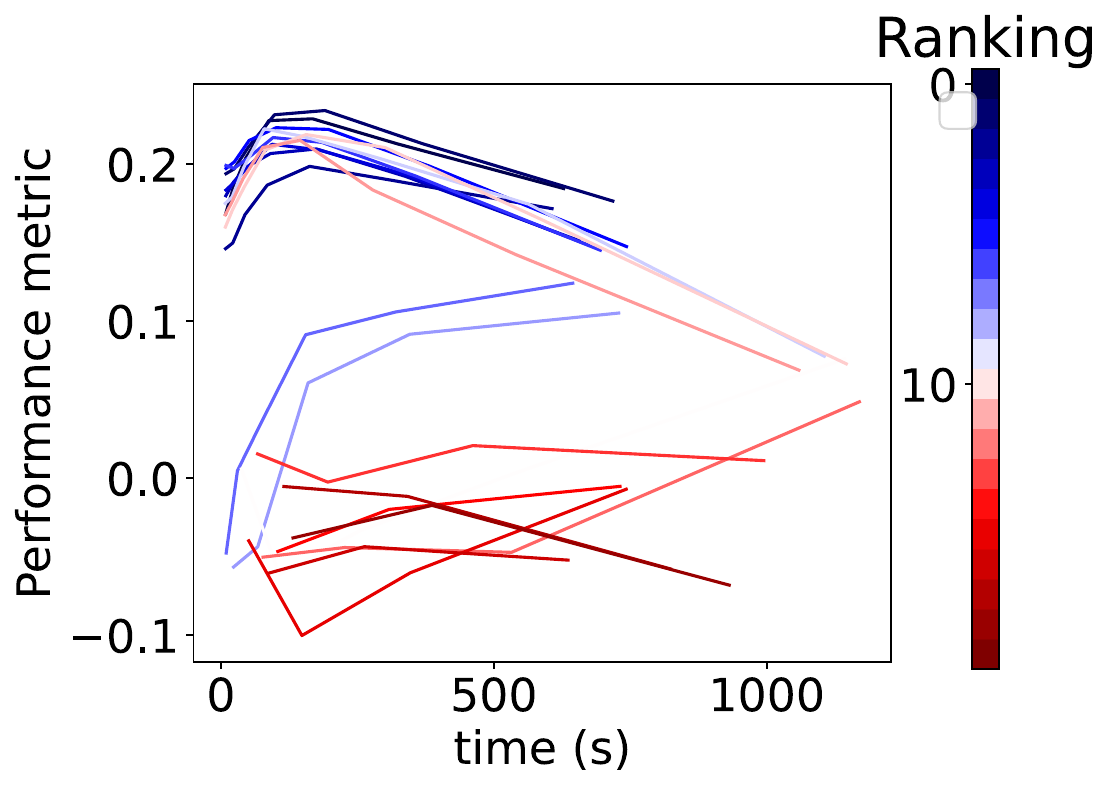}}
        \caption{{carlo}}
        
    \end{subfigure}
    \hfill
    \begin{subfigure}[b]{0.245\textwidth}
        \centering
        {\includegraphics[width=1\linewidth]{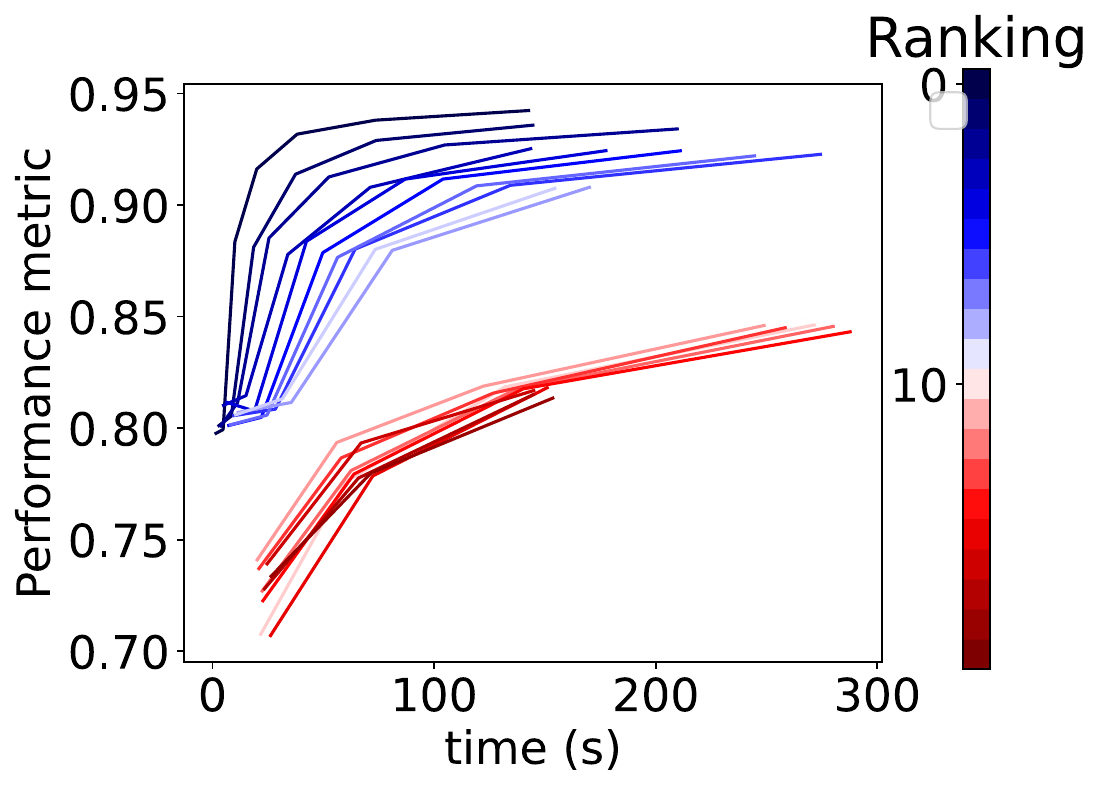}}
        \caption{{digits}}
       
    \end{subfigure}
    \hfill
    \begin{subfigure}[b]{0.245\textwidth}
        \centering
        {\includegraphics[width=1\linewidth]{dionis.pdf}}
        \caption{{dionis}}

    \end{subfigure}

    \begin{subfigure}[b]{0.245\textwidth}
        \centering
        {\includegraphics[width=1\linewidth]{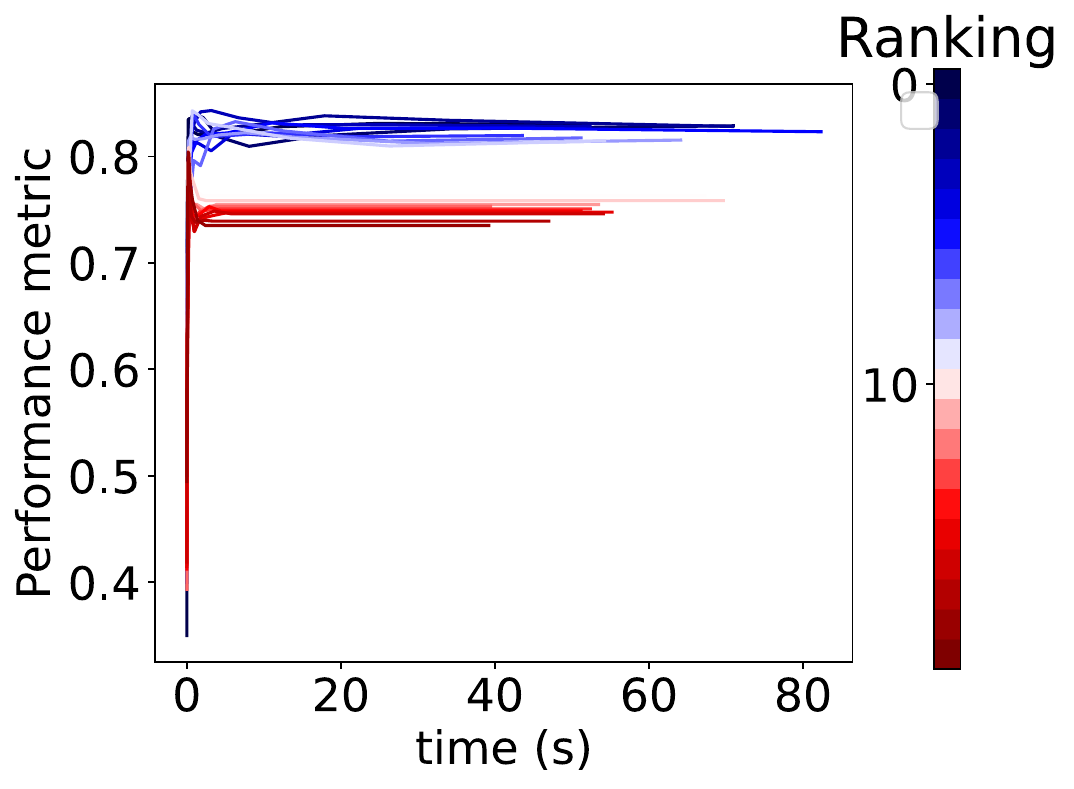}}
        \caption{{dorothea}}
        
    \end{subfigure}
    \hfill
    \begin{subfigure}[b]{0.245\textwidth}
        \centering
        {\includegraphics[width=1\linewidth]{evita.pdf}}
        \caption{{evita}}
       
    \end{subfigure}
    \hfill
    \begin{subfigure}[b]{0.245\textwidth}
        \centering
        {\includegraphics[width=1\linewidth]{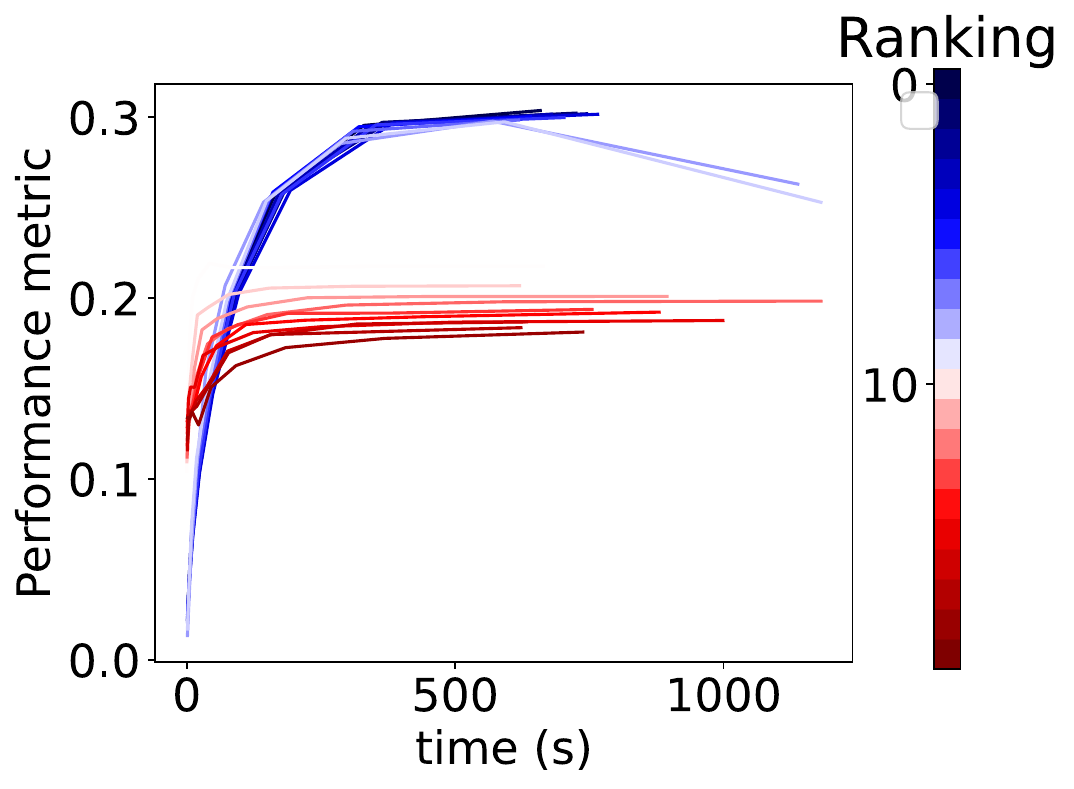}}
        \caption{{fabert}}
      
    \end{subfigure}
    \hfill
    \begin{subfigure}[b]{0.245\textwidth}
        \centering
        {\includegraphics[width=1\linewidth]{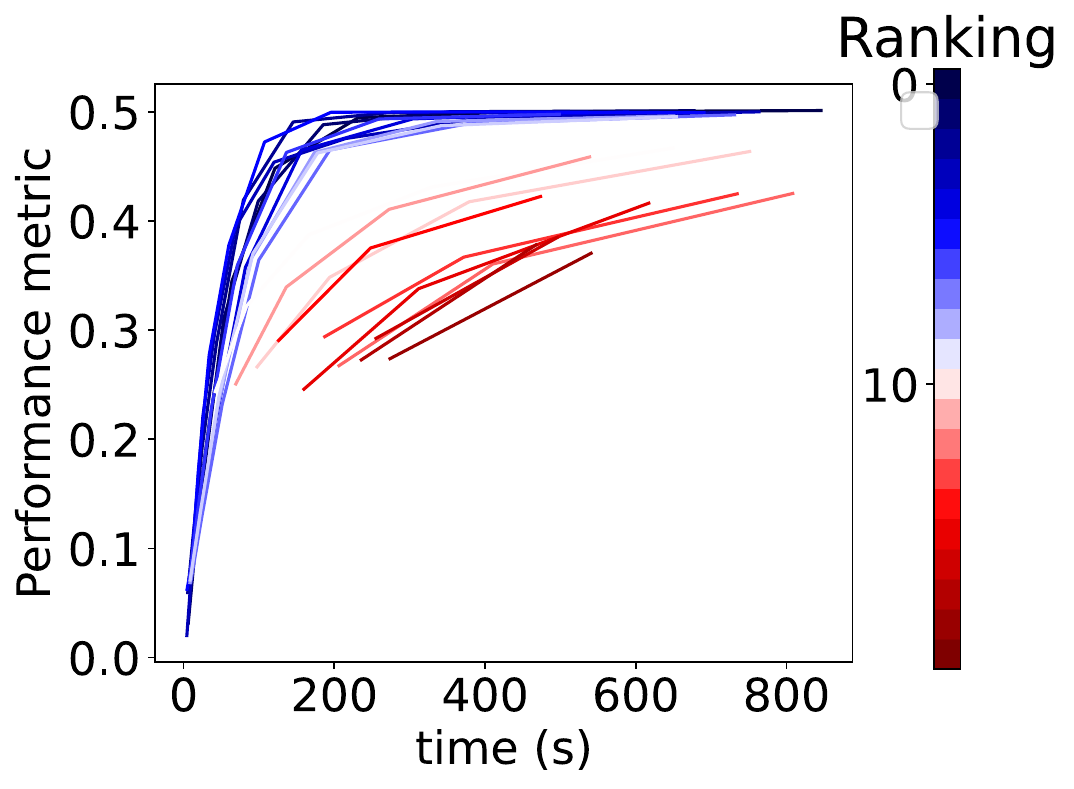}}
        \caption{{flora}}
       
    \end{subfigure}
    \caption{\textbf{Algorithms' learning curves with their final rankings}. We show some datasets where the baseline \textsc{BoS} beats the baseline \textsc{DDQN}. First round: (a-g); second round: (h-p). Algorithms are color-coded based on their final ranking (i.e., by comparing the last points on their learning curves). In these datasets, the learning curves do not cross each other very often, and the algorithm that ranked first early tends to maintain a very high rank at the end. This illustrates scenarios where \textsc{BoS} beats DDQN. However, in practice, one cannot know in advance if the learning curves of algorithms will cross each other often on a given dataset.} 
    \label{fig:algo-ranking}
\end{figure*}

\section{Conclusions}
{We addressed the machine learning algorithm selection problem under budget constraints. Our proposed framework, based on MDP, enables simultaneous algorithm selection and budget allocation by leveraging learning curves during the learning process. We organized challenges using novel learning curve datasets, allowing agents to meta-learn from knowledge gained on other datasets. The challenges' results demonstrated that agents using meta-learned knowledge of learning curve progressions outperformed those without this capability, demonstrated by their action trajectories.} The challenges will remain open for post-challenge submissions and serve as long-lasting benchmarks. Promising future directions include the meta-training of more sophisticated Reinforcement Learning (RL) methods, such as Proximal Policy Optimization by \cite{DBLP:journals/corr/SchulmanWDRK17}, and the expansion of our benchmarks to encompass a wider variety of algorithms and types of learning curves. {Another future research direction is to train multiple meta-learning agents to collaboratively select algorithms and budgets, motivated by the demonstrated efficacy of ensembles of algorithm selectors by \cite{DBLP:journals/ml/TornedeGTWH23}.}

\section*{Acknowledgements} The authors would like to thank Romain Egele and Felix Mohr for fruitful discussions, Nathan Grinsztajn and Adrien Pavao for helping organize the challenge, and Jan van Rijn for beta-testing it. We extend our appreciation to the challenge winners, who open-sourced their methods and provided feedback. This work was supported by ANR Chair of Artificial Intelligence HUMANIA ANR-19-CHIA-0022 and TAILOR EU Horizon 2020 grant 952215.

\bibliographystyle{abbrv}
\bibliography{refs}

\end{document}